\documentclass[final,12pt]{elsarticle}

\usepackage{amsmath}
\usepackage{amsthm}
\usepackage{amssymb}
\usepackage{epsfig}
\usepackage{mathrsfs}
\usepackage{graphicx}

\usepackage{caption}
\usepackage{subcaption}

\usepackage{abceqn}

\usepackage{fancyheadings}
\usepackage{tensor}
\usepackage{mathtools}
\usepackage{natbib}
\usepackage{ulem}
\usepackage{amssymb}
\usepackage{multirow}
\usepackage{hhline}
\usepackage{epstopdf}

\usepackage{color}
\usepackage{verbatim}

\usepackage{subcaption}

\newtheorem{thm}{Theorem}[section]

\theoremstyle{remark}
\newtheorem{rem}[thm]{Remark}

\theoremstyle{definition}

\usepackage{pslatex}

\newenvironment{mymathbox}
{\par\smallskip\centering\begin{lrbox}{0}%
\begin{minipage}[c]{1.05\textwidth}}
{\end{minipage}\end{lrbox}%
\framebox[1.13\textwidth]{\usebox{0}}%
\par\medskip
\ignorespacesafterend}

\newenvironment{mymathbox2}
{\par\smallskip\centering\begin{lrbox}{0}%
\begin{minipage}[c]{0.4\textwidth}}
{\end{minipage}\end{lrbox}%
\framebox[0.6\textwidth]{\usebox{0}}%
\par\medskip
\ignorespacesafterend}

\journal{ArXiv }

\begin{document}

\begin{frontmatter}
	
\title{A Petrov-Galerkin Spectral Element Method for Fractional Elliptic Problems}
\author{Ehsan Kharazmi$\,^{a,b}$, Mohsen Zayernouri$\,^{a,b,*}$ and George Em Karniadakis$\,^{c}$}

\address{$^a$ Department of Computational Mathematics, Science, and Engineering $\&$\\ 
		$^b$Department of Mechanical Engineering\\
		Michigan State University, 428 S Shaw Ln, East Lansing, MI 48824}

\address{$^c$ D\lowercase{ivision of} A\lowercase{pplied} M\lowercase{athematics}, B\lowercase{rown} U\lowercase{niversity}, 182 G\lowercase{eorge}, P\lowercase{rovidence}, RI 02912, USA}
	
\cortext[cor1]{Corresponding Author: zayern$@$msu.edu}
	
%\maketitle
%\tableofcontents
\begin{abstract}
We develop a new $C^{\,0}$-continuous Petrov-Galerkin spectral element method for one-dimensional fractional elliptic problems of the form $\prescript{}{0}{\mathcal{D}}_{x}^{\alpha} u(x) - \lambda u(x) = f(x)$, $\alpha \in (1,2]$, subject to homogeneous boundary conditions. We employ the standard (modal) spectral element bases and the Jacobi poly-fractonomials as the test functions \cite{Zayernouri2013}. We formulate a new procedure for assembling the global linear system from elemental (local) mass and stiffness matrices. The Petrov-Galerkin formulation requires performing elemental (local) construction of mass and stiffness matrices in the standard domain only once. Moreover, we efficiently obtain the non-local (history) stiffness matrices, in which the non-locality is presented analytically for uniform grids. We also investigate two distinct choices of basis/test functions: i) local basis/test functions, and ii) local basis with global test functions. We show that the former choice leads to a better-conditioned system and accuracy. We consider smooth and singular solutions, where the singularity can occur at boundary points as well as in the interior domain. We also construct two non-uniform grids over the whole computational domain in order to capture singular solutions. Finally, we perform a systematic numerical study of non-local effects via full and partial history fading in order to further enhance the efficiency of the scheme.
\end{abstract}
\end{frontmatter}

%\tableofcontents

%\newpage

%%%%%%%%%%%%%%%%%%%%%%%%%%%%%%%%%%%%%%%%%%%%%%%%%%%%%%%%%%%%%%%%%%%%%%%%%%%%%%%%%%%%%%%%%%%%%%%%%%%%%%%%%%%%%%%%%%%%%%%%%%%%%%%%%%%%%%%%%%%%%%%%%%%%%%%%%%%%
%%%%%%%%%%%%%%%%%%%%%%%%%%%%%%%%%%%%%%%%%%%%%%%%%%%%%%%%%%%%%%%%%%%%%%%%%%%%%%%%%%%%%%%%%%%%%%%%%%%%%%%%%%%%%%%%%%%%%%%%%%%%%%%%%%%%%%%%%%%%%%%%%%%%%%%%%%%%
%
%%%%%%%%%%%%%%%%%%%%%%%%%%%%%%%%%%
%
\section{Introduction}
\label{Sec: Introduction}
%
%%%%%%%%%%%%%%%%%%%%%%%%%%%%%%%%%%
%

Fractional order models open up new possibilities for robust mathematical modeling of complex multi-scale problems and anomalous transport phenomena including: non-Markovian (L\'{e}vy flights) processes in turbulent flows \cite{jha2003evidence, Castillo2004Plasma}, non-Newtonian fluids and rheology \cite{jaishankar2013}, non-Brownian transport phenomena in porous and disordered materials \cite{meer01, meral2010fractional}, non-Gaussian processes in multi-scale complex fluids and multi-phase applications \cite{jaishankar2014}, visco-elastic bio-tissues, and visco-elasto-plastic materials \cite{meral2010fractional,naghibolhosseini2015estimation,suzuki2016fractional}.

A number of local numerical methods, prominently finite difference methods (FDMs), have been developed for solving fractional partial differential equations (FPDEs) \cite{Lubich1983, Lubich1986, Sanz1988, Sugimoto1991, Metzler2000, Gorenflo2002, Diethelm2004, Langlands2005, Sun2006,  Lin2007, wang2010direct, wang2011fast, Huang2012, Cao2013, zeng2015numerical,zayernouri2016fractional}. Fix and Roop \cite{Fix2004} developed the first theoretical framework for the least-square finite element method (FEM) approximation of a fractional-order differential equation, where optimal error estimates are proven for piecewise linear elements. However, Roop \cite{Roop2006} later showed that the main hurdle to overcome in the FEM is the non-local nature of the fractional operator, which leads to large dense matrices; he showed that even the construction of such matrices presents difficulties. Ervin and Roop \cite{ErvRoop06} presented a theoretical framework for the variational solution of the steady state fractional advection dispersion equation based on FEM and proved the existence and uniqueness of the results. Jin et al. \cite{jin2014error} proved the existence and uniqueness of a weak solution to the space-fractional parabolic equation using FEM; they showed an enhanced regularity of the solution and derived the error estimate for both semidiscrete and fully discrete solution. Well-posedness, regularity of the weak solution, stability of the discrete variational formulation and error estimate of the FEM approximation were investigated for fractional elliptic problems in \cite{jin2015variational}. Wang and Yang \cite{wang2013wellposedness} generalized the analysis to the case of fractional elliptic problems with variable coefficient, analyzed the regularity of the solution in H$\ddot{\text{o}}$lder spaces, and established the well-posedness of a Petrov-Galerkin formulation. Wang et al. \cite{wang2014inhomogeneous} developed an indirect FEM for the Dirichlet boundary-value problems of Caputo FPDEs showing the reduction in the computational work for numerical solution and memory requirements.

There has been recently more attention and effort put on developing global and high-order approximations, which are capable of efficiently capturing the inherent non-local effects. A Chebyshev spectral element method (SEM) for fractional-order transport was adopted by Hanert \cite{Hanert2010} and later on, the idea of least-square FEM was extended to SEM by Carella \cite{Carella2012}. More recently, Deng and Hesthevan \cite{Deng2013} and Xu and Hesthaven \cite{xu2014discontinuous} developed local DG methods for solving space-fractional diffusion and convection-diffusion problems.

Two new spectral theories on fractional and tempered fractional Sturm-Liouville problems (TFSLPs) have been developed by Zayernouri et al. in \cite{Zayernouri2013, zayernouri2015tempered}. This approach first fractionalizes and then tempers the well-known theory of Sturm-Liouville eigen-problems. The explicit eigenfunctions of TFSLPs are analytically obtained in terms of \textit{tempered Jacobi poly-fractonomials}. These poly-fractonomials have been successfully employed in developing a series of high-order and efficient Petrov-Galerkin spectral and discontinuous spectral element methods \cite{Zayernouri14-SIAM-Frac-Advection,Zayernouri2015Unified,ZayernouriVariableOrder2015JCP}. In \cite{Zayernouri_FODEs_2014}, Zayernouri and Karniadakis developed a spectral and spectral element method for FODEs with an exponential accuracy. They also developed a highly accurate discontinuous SEM for time- and space- fractional advection equation in \cite{Zayernouri14-SIAM-Frac-Advection}. Dehghan et al. \cite{dehghan2016legendre} considered Legendre SEM in space and FDM in time for solving time-fractional sub-diffusion equation. Su \cite{su2015parallel} provided a parallel spectral element method for the fractional Lorenz system and a comparison of the method with FEM and FDM.

%\cite{Zayernouri2015_FracAdams_Chemotaxis}

The SEM discretization has the benefit of domain decomposition into non-overlapping elements, which potentially provide a geometrical flexibility, especially for adaptivity as well as complex domains. Moreover, high-order approximations within each element yield a fast rate of convergence even in the cases of non-smooth and/or rapid transients in the solution. Therefore, a tractable computational cost of the method can be achieved by a successful combination of \textit{h-refinement}, where the solution is rough, and \textit{p-refinement}, where the solution is smooth.

In the present work, we consider the one-dimensional space-fractional Helmholtz equation of order $\alpha \in (1,2]$  subject to homogeneous boundary conditions. We formulate a weak form, in which the fractional portion $\mu \in (0,1]$ is transfered onto some proper fractional order test functions via integration-by-parts. This setting enables us to employ the standard polynomial modal basis functions, used in SEM \cite{Karniadakis2005}. Subsequently, we develop a new  $C^{\,0}$-continuous Petrov-Galerkin SEM, following the recent spectral theory of fractional Sturm-Liouville problem, where the test functions are of Jacobi poly-fractonomials of second kind \cite{Zayernouri2013}. We investigate two distinct choices of basis/test functions: i) \textit{local} basis/test functions, and ii) \textit{local} basis with \textit{global} test functions, which enables the construction of elemental mass/stiffness matrices in the standard domain $[-1,1]$. We explicitly compute the elemental stiffness matrices using the orthogonality of Jacobi polynomials. Moreover, we efficiently obtain the non-local (history) stiffness matrices, in which the non-locality is presented \textit{analytically}. On one hand, we formulate a new \textit{non-local assembling procedure} in order to construct the global linear system from the local (elemental) mass/stiffness matrices and history matrices. On the other hand, we formulate a procedure for \textit{non-local scattering} to obtain the elemental expansion coefficients from the global degrees of freedom. We demonstrate the efficiency of the Petrov-Galerkin methods and show that the choice of local bases/test functions leads to a better accuracy and conditioning. Moreover, for uniform grids, we compute the history matrices off-line. The stored history matrices can be retrieved later in the construction of the global linear system. We show the great improvement in the computational cost by performing the retrieval procedure compared to on-line computation. We also introduce a non-uniform \textit{kernel-based} grid generation in addition to \textit{geometrically progressive} grid generation approaches. Furthermore, we investigate the performance of the developed schemes by considering two cases of smooth and singular solutions, where the singularity can occur at boundary points or the interior domain. Finally, we study the effect of history fading via a systematic analysis, where we consider the history up to some specific element and let the rest fade. This results in less computational cost, while we show that the accuracy is still preserved. The main contributions of this work are listed in the following:
\begin{itemize}
	\item Development of a new fast and accurate $C^{\,0}$-continuous Petrov-Galerkin spectral element method, employing local basis/test functions, where the test functions are Jacobi poly-fractonomials.
	
	\item Reducing the number of history calculation from $\frac{N_{el} (N_{el}-1)}{2}$ to $(N_{el}-1)$ for a uniformly partitioned domain.
	
	\item Analytical expression of non-local effects in uniform grids leading to fast computation of the history matrices.
	
	\item A new procedure for the assembly of the global linear system.
	
	\item Performing off-line computation of history matrices and on-line retrieval of  the stored matrices.
	
	\item Boundary and interior singularity capturing using adaptive \textit{hp}-refinement.
	
	\item Non-uniform kernel-based grid generation for resolving steep gradients and singularities.
\end{itemize}

%further investigate the efficiency of the method through numerical examples for different number of elements and modes, where we show that in the case of local test functions, the assembled linear system has a lower condition number and thus leads to a better approximability. Moreover, to study the effect of non-locality, we compute the $L_2$-norm error of approximation by considering a systematic history fading. in the \textit{fully} history fading, the effect of non-locality is taken into account up to some past element and then the rest of history is truncated; whereas, in the \textit{partially} history fading the rest of history is partially computed and not fully truncated.

The organization of the paper is as follows: section \ref{Sec: Definitions} provides preliminary definitions including problem definition, derivation of the weak form and expressions for the local basis and local/global test functions. In section \ref{Sec: PG Local Test}, we present a Petrov-Galerkin method employing the local basis/test functions in addition to formulating the non-local assembling and non-local scattering procedures, followed by a discussion on how to compute \textit{off-line} the history matrices. We also present the two non-uniform grid generation approaches. In section \ref{Sec: PG Global Test}, we present a Petrov-Galerkin method employing the local basis with global test functions, compared with the former scheme. In section \ref{Sec: NumExamp}, we demonstrate the computational efficiency of the methods by considering several numerical examples of smooth and singular solutions. Finally, we perform the off-line computation and retrieval procedure of history matrices and a systematic history fading analysis. We end the paper with a summary.

%%%%%%%%%%%%%%%%%%%%%%%%%%%%%%%%%%%%%%%%%%%%%%%%%%%%%%%%%%%%%%%%%%%%%%%%%%%%%%%%%%%%%%%%%%%%%%%%%%%%%%%%%%%%%%%%%%%%%%%%%%%%%%%%%%%%%%%%%%%%%%%%%%%%%%%%%%%%
%%%%%%%%%%%%%%%%%%%%%%%%%%%%%%%%%%%%%%%%%%%%%%%%%%%%%%%%%%%%%%%%%%%%%%%%%%%%%%%%%%%%%%%%%%%%%%%%%%%%%%%%%%%%%%%%%%%%%%%%%%%%%%%%%%%%%%%%%%%%%%%%%%%%%%%%%%%%

%
%%%%%%%%%%%%%%%%%%%%%%%%%%%%%%%%%%
%
\section{Definitions}
\label{Sec: Definitions}
%
%%%%%%%%%%%%%%%%%%%%%%%%%%%%%%%%%%
%

Let $ \xi \in [-1,1]$. Then, the left-sided and right-sided Riemann-Liouville integrals of order $\sigma$,  $n-1 < \sigma \leq n$, $n \in \mathbb{N}$, are defined (see e.g., \cite{Miller93, Podlubny99}) respectively as
\begin{equation}
	\label{Eq: left RL integral}
	(\prescript{RL}{-1}{\mathcal{I}}_{\xi}^{\sigma}) u(\xi) = \frac{1}{\Gamma(\sigma)} \int_{-1}^{\xi} \frac{u(s) ds}{(\xi - s)^{n-\sigma} },\quad \xi>-1 ,
\end{equation}
and
\begin{equation}
	\label{Eq: right RL integral}
	(\prescript{RL}{\xi}{\mathcal{I}}_{1}^{\sigma}) u(\xi) = \frac{1}{\Gamma(\sigma)} \int_{\xi}^{1} \frac{u(s) ds}{(s - \xi)^{n-\sigma} },\quad \xi<1.
\end{equation}
The corresponding left-sided and right-sided fractional derivatives of order $\sigma$ are then defined as  
\begin{equation}
	\label{Eq: left RL derivative}
	(\prescript{RL}{-1}{\mathcal{D}}_{\xi}^{\sigma}) u(\xi) = \frac{d^n}{d\xi^n} (\prescript{RL}{-1}{\mathcal{I}}_{\xi}^{n-\sigma} u) (\xi) = \frac{1}{\Gamma(n-\sigma)}  \frac{d^{n}}{d\xi^n} \int_{-1}^{\xi} \frac{u(s) ds}{(\xi - s)^{\sigma+1-n} },\quad \xi >-1 ,
\end{equation}
and
\begin{equation}
	\label{Eq: right RL derivative}
	(\prescript{RL}{\xi}{\mathcal{D}}_{1}^{\sigma}) u(\xi) = \frac{(-d)^n}{d\xi^n} (\prescript{RL}{\xi}{\mathcal{I}}_{1}^{n-\sigma} u) (\xi) = \frac{1}{\Gamma(n-\sigma)}  \frac{(-d)^{n}}{d\xi^n} \int_{\xi}^{1} \frac{u(s) ds}{(s - \xi)^{\sigma+1-n} },\quad \xi < 1 ,
\end{equation}
respectively.

By performing an affine mapping from the standard domain $[-1,1]$ to the interval $x \in [x_{\varepsilon - 1},x_{\varepsilon}]$, we obtain
\begin{eqnarray}
	\label{Eq: RL in xL-xR}
	\prescript{RL}{x}{\mathcal{D}}_{x_{\varepsilon}}^{\mu} u  &=&  (\frac{2}{x_{\varepsilon} - x_{\varepsilon - 1}})^\mu (\prescript{RL}{-1}{\mathcal{D}}_{\xi}^{\mu} \, u )(\xi). 
	%\\ 
	%\label{Eq: Caputo in xL-xR}
	%\prescript{C}{a}{\mathcal{D}}_{x}^{\sigma} u  &=&  (\frac{2}{b-a})^\sigma (\prescript{C}{-1}{\mathcal{D}}_{\xi}^{\sigma} \, u) (\xi).
	%
\end{eqnarray} 
Hence, we can perform the operations in the standard domain only once for any given $\sigma$ and efficiently utilize them on any arbitrary interval without resorting to repeating the calculations.

We define the \textit{Jacobi poly-fractonomial}s (of second kind), used as the test functions in developing the proposed numerical schemes following the recent theory of fractional Sturm-Liouville eigen-problems (FSLP) in \cite{Zayernouri2013}. The corresponding regular poly-fractonomials are given in the standard domain $[-1,1]$ by
\begin{equation}
	\label{Eq: Jacobipoly II}
	\prescript{(2)}{}{ \mathcal{P}}_{k}^{\,\,\mu}(\xi) = (1-\xi)^{\mu} P_{k-1}^{\mu,-\mu} (\xi),\quad \xi \in [-1,1],
\end{equation}
where $P_{k-1}^{\mu,-\mu}$ is the Jacobi Polynomial.

%
%%%%%%%%%%%%%%%%%%%%%%%%%%%%%%%%%%%%%%%%%
\subsection{\textbf{Problem Definition}}
\label{Sec: Model Problem}
%%%%%%%%%%%%%%%%%%%%%%%%%%%%%%%%%%%%%%%%%
%
We study the following fractional Helmholtz equation of order $\alpha = 1+\mu$, $\mu \in (0,1]$:
\begin{eqnarray}
	\label{Eq: Helmholtz}
	\prescript{RL}{0}{\mathcal{D}}_{x}^{\alpha} u(x) - \lambda u(x) &=& f(x), \quad \forall x \in \Omega
	\\ 
	u(0)&=&u(L)= 0, \quad \forall x \in \partial \Omega,
\end{eqnarray}
where $\Omega = [0,L]$. By multiplying both sides of \eqref{Eq: Helmholtz} by some proper test function $v(x)$, then taking the fractional integration-by-parts, we obtain the following bilinear form:
\begin{eqnarray}
	\label{Eq: bilnear-from}
	a(u, v) = l(v),
\end{eqnarray}
in which
\begin{eqnarray}
	\label{Eq: a def}
	a(u, v) &=& \Big(\frac{du}{dx} \,,\, \prescript{RL}{x}{\mathcal{D}}_{L}^{\mu} \,v \Big)_{\Omega}- \lambda \Big( u \,,\, v\Big)_{\Omega}, 
	\\ 
	\label{Eq: l def}
	l(v) &=& \Big( f \,,\,v \Big)_{\Omega}, 
\end{eqnarray}
where $( \cdot \,,\,\cdot )_{\Omega}$ denotes the well-known inner-product. 
%
%%%%%%%%%%%%%%%%%%%%%%%%%%%%%%%%%%%%%%%%%
\subsection{\textbf{ Local Basis Functions}  }
\label{Subsec: Basis}
\begin{figure}[t]
	\center
	\includegraphics[width=1\linewidth]{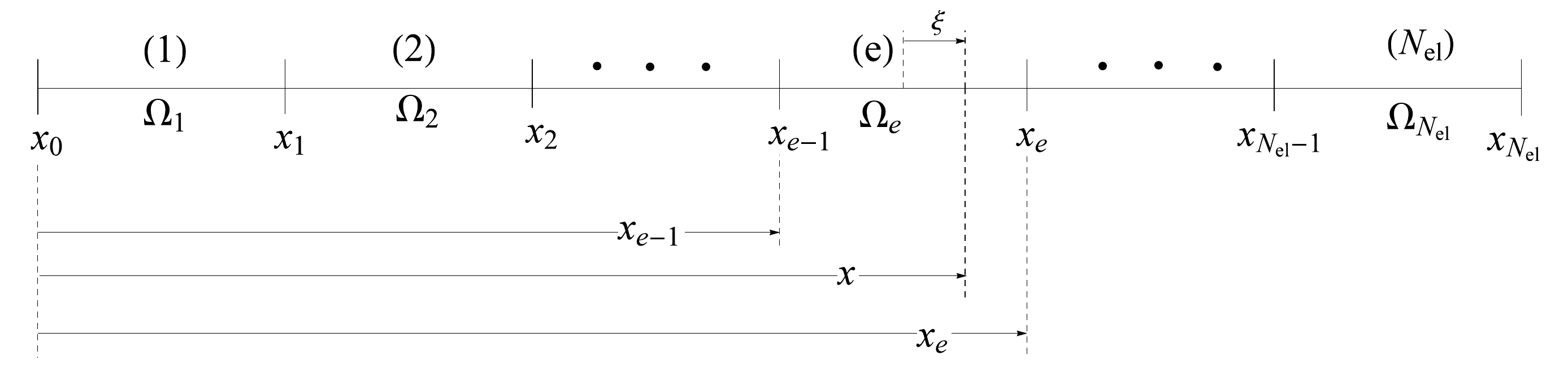}
	\caption{\label{Fig: Domain Disc} Domain partitioning}
\end{figure}

%%%%%%%%%%%%%%%%%%%%%%%%%%%%%%%%%%%%%%%%%
%
We partition the computational domain into $N_{el}$ non-overlapping elements $\Omega_e = \left[ x_{e-1} , x_{e} \right]$ such that $\Omega = \cup_{e=1}^{Nel} \Omega_e$, see Fig. \ref{Fig: Domain Disc}. 
%
%******************************************************************************************
Therefore, the bilinear form \eqref{Eq: a def} can be written as
\begin{align}
	\label{Eq: a def disc}
	a(u, v) \approx a(u^{\delta}, v^{\delta}) = \sum_{e=1}^{N_{el}} \Big(\frac{du_N^{(e)}}{dx} \,,\, \prescript{RL}{x}{\mathcal{D}}_{L}^{\mu} \,v^{\delta} \Big)_{\Omega_e}
	- \lambda  \sum_{e=1}^{N_{el}} \Big( u_N^{(e)} \,,\, v^{\delta} \Big)_{\Omega_e} ,
\end{align}
where we approximate the solution in each element as 
\begin{align}
	\label{Eq: u element-expansion}
	u_N^{(e)}(x) = \sum_{p=0}^{P} \hat{u}^{(e)}_p \psi_p (x), \quad x \in \Omega_{e},
\end{align}
and thus, the approximated solution over the whole domain is
\begin{eqnarray}
	\label{Eq: u-expansion}
	u \approx u^{\delta} (x) = \sum_{e=1}^{N_{el}} \sum_{p=0}^{P} \hat{u}^{(e)}_p \psi_p (x).
\end{eqnarray}
We choose the $P+1$ modal basis functions $\psi_p ( x )$, defined in the standard element in terms of $\zeta \in [-1,1]$ as
\begin{eqnarray}
	\label{Eq: basis}
	\psi_p ( \zeta  ) =
	\begin{cases}
		\frac{1- \zeta}{2},\quad & p=0, \\
		%\\
		(\frac{1- \zeta}{2})(\frac{1+ \zeta}{2}) P^{1,1}_{p-1}(\zeta),\quad & p=1,2,\cdots, P-1, \\
		%\\
		\frac{1+ \zeta}{2},\quad & p=P, 
	\end{cases}
\end{eqnarray}
which are also used in standard spectral element methods for integer-order PDEs (see e.g., \cite{Karniadakis2005}).
%
%%%%%%%%%%%%%%%%%%%%%%%%%%%%%%%%%%%%%%%%%
\subsection{ \textbf{Test Functions: Local vs. Global} }
\label{Subsec: Test}
%%%%%%%%%%%%%%%%%%%%%%%%%%%%%%%%%%%%%%%%%
%
We choose two types of test functions $v^{\delta}$: i) \textit{local} test functions, and ii) \textit{global} test functions, defined for $\varepsilon =1,2, \cdots, N_{el}$ as follows:
\begin{eqnarray}
	\label{Eq: test}
	v_k^{\,local} (x) = v_k^{\varepsilon} (x) =
	\begin{cases}
		\prescript{(2)}{}{\mathcal{P}}_{k+1}^{\mu}( x^{\varepsilon}  ) ,\quad & \forall x \in \Omega_{\varepsilon}, \\
		%\\
		0,\quad & otherwise,
	\end{cases}
	,\quad k=0,1,\cdots,P ,
\end{eqnarray}
in which $\prescript{(2)}{}{\mathcal{P}}_{k+1}^{\mu}( x^{\varepsilon}  )$ represents the Jacobi \textit{poly-fractonomial} of second kind, defined in the corresponding intervals $\Omega_{\varepsilon} = [x_{\varepsilon -1} , x_{\varepsilon}]$, using \eqref{Eq: Jacobipoly II}; and
\begin{eqnarray}
	\label{Eq: global test}
	v_k^{\,global} (x) = v_k^{\varepsilon} (x) =
	\begin{cases}
		\prescript{(2)}{}{\mathcal{P}}_{k+1}^{\mu}( x^{1\sim  \varepsilon}  ) ,\quad & \forall x \in [0 , x_{\varepsilon}],
		\\
		0,\quad & otherwise
	\end{cases}
	,\quad k=0,1,\cdots,P ,
\end{eqnarray}
where $\prescript{(2)}{}{\mathcal{P}}_{k+1}^{\mu}( x^{1\sim  \varepsilon}  ) $ represents the Jacobi \textit{poly-fractonomial} of second kind, defined in the corresponding intervals $[0 , x_{\varepsilon}]$, using \eqref{Eq: Jacobipoly II}. It should be noted that for each element $\varepsilon$, the corresponding local test function has nonzero value only in the element and vanishes elsewhere, unlike the corresponding global test function, which vanishes only where $x > x_{\varepsilon}$.

%%%%%%%%%%%%%%%%%%%%%%%%%%%%%%%%%%%%%%%%%%%%%%%%%%%%%%%%%%%%%%%%%%%%%%%%%%%%%%%%%%%%%%%%%%%%%%%%%%%%%%%%%%%%%%%%%%%%%%%%%%%%%%%%%%%%%%%%%%%%%%%%%%%%%%%%%%%%%%%%%%%%%%%%%%%%%%%%%%%%%%%%%%%%%%%%%%%%%%%%%%%%%%%%%%%%%%%%%%%%%%%%%%%%%%%%%%%%%%%%%%%%%%%%%%%%%%%%%%%%%%%%%%%%%%%%%%%%%%%%%%%%%%%%%%%%%%%%%%%%%%%%%%%%%%%%

%
%%%%%%%%%%%%%%%%%%%%%%%%%%%%%%%%%%%%%%%%%
\section{Petrov-Galerkin Method with Local Test Functions}
\label{Sec: PG Local Test}
\begin{figure}[t]
	\center
	\includegraphics[width=0.65\linewidth]{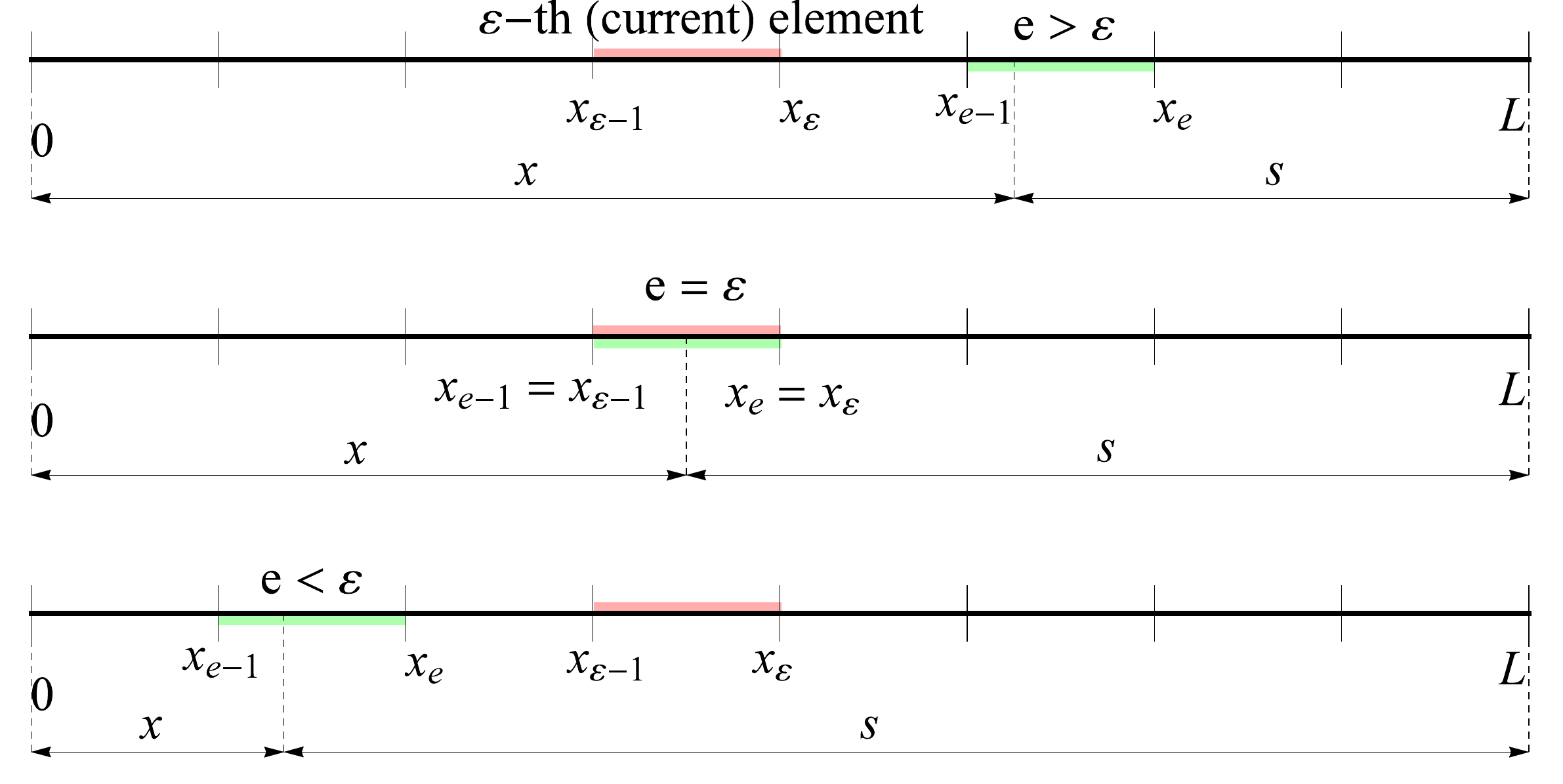}
	\caption{\label{Fig: history function} Location of the (dummy) element number, $e$, with respect to the current element, $\varepsilon$. If $e > \varepsilon$, (top), then $\prescript{RL}{x}{\mathcal{D}}_{L}^{\mu} \,v_k^{\varepsilon} (x) = 0$. If $e = \varepsilon$, (middle), then $\prescript{RL}{x}{\mathcal{D}}_{L}^{\mu} \,v_k^{\varepsilon} (x) =  \prescript{RL}{x}{\mathcal{D}}_{x_{\varepsilon}}^{\mu} \Big[ \prescript{(2)}{}{\mathcal{P}}_{k+1}^{\mu}(x) \, \Big]$. If $e < \varepsilon$, (bottom), then $\prescript{RL}{x}{\mathcal{D}}_{L}^{\mu} \,v_k^{\varepsilon} (x) = H_k^{(\varepsilon)}(x)$. }
\end{figure}

%%%%%%%%%%%%%%%%%%%%%%%%%%%%%%%%%%%%%%%%%
%
By substituting \eqref{Eq: u element-expansion} and \eqref{Eq: test} into \eqref{Eq: a def disc}, we obtain:
\begin{eqnarray}
	\nonumber
	&&
	\sum_{e=1}^{N_{el}}  
	\Big(  \sum_{p=0}^{P} \hat{u}^{(e)}_p \,\, \frac{d\psi_p(x) }{dx} \,,\, \prescript{RL}{x}{\mathcal{D}}_{L}^{\mu} \,v_k^{\varepsilon} (x) \Big)_{\Omega_e}  
	-
	\lambda 
	\sum_{e=1}^{N_{el}} 
	\Big( 
	\sum_{p=0}^{P} \hat{u}^{(e)}_p \,\, \psi_p (x) \,,\, \,v_k^{\varepsilon} (x)
	\Big)_{\Omega_e}  
	\\ \label{Eq: a def disc subs-1}
	&=&   \Big( f \,,\,v_k^{\varepsilon} (x) \Big)_{\Omega}, \quad \varepsilon=1,2,\cdots,N_{el}, \quad k=0,1,\cdots,P .
\end{eqnarray}
Since the local test function vanishes $\forall x \in \Omega_e \neq \Omega_{\varepsilon}$, we have 
\begin{eqnarray}
	\nonumber
	\lambda
	\sum_{e=1}^{N_{el}}
	\Big( \sum_{p=0}^{P} \hat{u}^{(e)}_p \psi_p ( x  ) \,,\, \,v_k^{\varepsilon} (x) \Big)_{\Omega_e}  
	&= &  
	\lambda
	\Big(\sum_{p=0}^{P} \hat{u}^{(\varepsilon)}_p \psi_p ( x  ) \,,\, \,v_k^{\varepsilon} (x) \Big)_{\Omega_{\varepsilon}}, 
	\\ \nonumber
	\Big( f \,,\,v_k^{\varepsilon} (x) \Big)_{\Omega} &= &
	\Big( f \,,\,v_k^{\varepsilon} (x) \Big)_{\Omega_{\varepsilon}}.
\end{eqnarray}
Moreover, for every $\varepsilon$, the right-sided fractional derivative, 
\begin{align*}
	\prescript{RL}{x}{\mathcal{D}}_{L}^{\mu} \,v_k^{\varepsilon} (x) 
	= \frac{-1}{\Gamma(1-\mu)} \frac{d}{dx} \int_{x}^{L}  \frac{v_k^{\varepsilon} (s) }{(s - x)^{\mu}} ds, \quad x \in \Omega_e,
\end{align*}
is taken from $x \in \Omega_e$ to $x=L$, where $e = 1,2,\cdots, N_{el}$ through the summation over the elements and $s$ varies from $x \in \Omega_e$ to $L$. The local test function vanishes $\forall x \in \Omega_e \neq \Omega_{\varepsilon}$, thus if $e > \varepsilon$ ($x > x_{\varepsilon}$, see Fig. \ref{Fig: history function} top), then
\begin{align}
	\label{Eq: frac der test function-1}
	\prescript{RL}{x}{\mathcal{D}}_{L}^{\mu} \,v_k^{\varepsilon} (x) 
	= \frac{-1}{\Gamma(1-\mu)} \frac{d}{dx} \int_{x}^{L}  \frac{ 0 }{(s - x)^{\mu}} ds
	= 0 ,
\end{align}
and if $e < \varepsilon$ ($x < x_{\varepsilon-1}$, see Fig. \ref{Fig: history function} bottom), then
\begin{align}
	\label{Eq: frac der test function-2}
	\prescript{RL}{x}{\mathcal{D}}_{L}^{\mu} \,v_k^{\varepsilon} (x) 
	= \frac{-1}{\Gamma(1-\mu)} \frac{d}{dx} \int_{x_{\varepsilon -1}}^{x_{\varepsilon}}  \frac{\prescript{(2)}{}{\mathcal{P}}_{k+1}^{\mu}(s)  }{(s - x)^{\mu}} ds
	\equiv H_k^{(\varepsilon)}(x) ,
\end{align}
and if $e = \varepsilon$, ($ x_{\varepsilon-1} < x < x_{\varepsilon}$, see Fig. \ref{Fig: history function} middle), then
\begin{align}
	\label{Eq: frac der test function-3}
	\prescript{RL}{x}{\mathcal{D}}_{L}^{\mu} \,v_k^{\varepsilon} (x) 
	= \frac{-1}{\Gamma(1-\mu)} \frac{d}{dx} \int_{x}^{x_{\varepsilon}}  \frac{\prescript{(2)}{}{\mathcal{P}}_{k+1}^{\mu}(s)  }{(s - x)^{\mu}} ds
	= \prescript{RL}{x}{\mathcal{D}}_{x_{\varepsilon}}^{\mu} \Big[ \prescript{(2)}{}{\mathcal{P}}_{k+1}^{\mu}(x) \, \Big] .
\end{align}
%
%
%******************************************************************************************

\noindent Hence, for $\varepsilon =1,2, \cdots, N_{el}$ and $k = 0 , 1 , \cdots , P$,
\begin{eqnarray}
	\label{Eq: frac der test function-4}
	\prescript{RL}{x}{\mathcal{D}}_{L}^{\mu} \,v_k^{\varepsilon} (x) =
	\begin{cases}
		0, \quad  & \forall x \in \Omega_e, \,\, e > \varepsilon,
		\\ 
		\\
		\prescript{RL}{x}{\mathcal{D}}_{x_{\varepsilon}}^{\mu} \,
		\Big[ \prescript{(2)}{}{\mathcal{P}}_{k+1}^{\mu}(x) \, \Big]  , \quad & \forall x \in \Omega_e, \,\, e = \varepsilon,
		% \forall x \in  [x_{\varepsilon -1} , x_{\varepsilon}]= \Omega_{\varepsilon},
		%
		\\ 
		\\
		H_k^{(\varepsilon)}(x) ,
		%\equiv
		%\frac{-1}{\Gamma(1-\mu)} \frac{d}{dx} \int_{x_{\varepsilon -1}}^{x^{\varepsilon}}  \frac{\prescript{(2)}{}{\mathcal{P}}_{k+1}^{\mu}(s)  }{(s - x)^{\mu}} ds,
		\quad & \forall x \in \Omega_e, \,\, e < \varepsilon.
		%%\forall x< x_{\varepsilon-1}.
	\end{cases}
	\vspace{0.5 in}
\end{eqnarray}

\noindent Therefore, the bilinear form \eqref{Eq: a def disc subs-1} can be written as
\begin{align}
	\nonumber
	&
	\sum_{e=1}^{\varepsilon - 1}  
	\sum_{p=0}^{P} \hat{u}^{(e)}_p \,\, 
	\Big(  \frac{d\psi_p(x) }{dx} \,,\, H_k^{(\varepsilon)}(x) \Big)_{\Omega_e}  
	+
	\sum_{p=0}^{P} \hat{u}^{(\varepsilon)}_p \,\, 
	\Big(  \frac{d\psi_p(x) }{dx} \,,\, \prescript{RL}{x}{\mathcal{D}}_{x_{\varepsilon}}^{\mu} \,
	\Big[ \prescript{(2)}{}{\mathcal{P}}_{k+1}^{\mu}(x) \, \Big] \Big)_{\Omega_\varepsilon} 
	\\
	\label{Eq: a def disc subs-2}
	-
	& \lambda \sum_{p=0}^{P} \hat{u}^{(\varepsilon)}_p 
	\Big(\psi_p ( x  ) \,,\, \, \prescript{(2)}{}{\mathcal{P}}_{k+1}^{\mu}(x) \Big)_{\Omega_{\varepsilon}}
	=   \Big( f \,,\, \prescript{(2)}{}{\mathcal{P}}_{k+1}^{\mu}(x) \Big)_{\Omega_{\varepsilon}},
\end{align}
and the weak form is obtained as
%%
%\begin{eqnarray}
%\nonumber
%%
%\sum_{e=1}^{N_{\varepsilon -1}} \sum_{p=0}^{P} \hat{u}^{(e)}_p 
%%
%\Big(\frac{d\psi^{e}_p }{dx} \,,\, H^{(\varepsilon)}_k(x)   \Big)_{\Omega_e} 
%&+&
%%
%(\frac{2}{x_{\varepsilon} - x_{\varepsilon} })^{\mu} \frac{\Gamma[1+k+\mu]}{\Gamma[1+k]} 
%%
%\sum_{p=0}^{P} \hat{u}^{(\varepsilon)}_p 
%\Big(\frac{d\psi^{\varepsilon}_p }{dx} \,,\, P_k(x^{\varepsilon})  \Big)_{\Omega_{\varepsilon}} 
%%
%\\ \nonumber
% &-&
%% 
% \lambda 
%\sum_{p=0}^{P} \hat{u}^{(\varepsilon)}_p  
% \Big( \psi^{\varepsilon}_p ( x  ) \,,\, \,\prescript{(2)}{}{\mathcal{P}}_{k+1}^{\mu}( x^{\varepsilon} \Big)_{\Omega_{\varepsilon}}  
%%
%%
%\\ \nonumber
%&=&  \Big( f \,,\,\prescript{(2)}{}{\mathcal{P}}_{k+1}^{\mu}( x^{\varepsilon}  ) \Big)_{\Omega_{\varepsilon}}, 
%%
%\end{eqnarray}
%%
%or in other words:
%
\begin{eqnarray}
	\label{Eq: Matrix Form}
	\sum_{e=1}^{\varepsilon -1} \sum_{p=0}^{P} \hat{u}^{(e)}_p 
	\,\,
	\hat{\textbf{S}}^{(e,\varepsilon)}_{kp}
	+
	\sum_{p=0}^{P} \hat{u}^{(\varepsilon)}_p 
	\Big[
	\,\,
	\textbf{S}^{(\varepsilon)}_{kp}
	- 
	\lambda 
	\,\,
	\textbf{M}^{(\varepsilon)}_{kp}
	\Big] = \textbf{f}^{(\varepsilon)}_k ,
	\quad
	\left\lbrace
	\begin{aligned}
		& \varepsilon=1,2,\cdots,N_{el}, 
		\\
		& k=0,1,\cdots,P ,
	\end{aligned}
	\right.
\end{eqnarray}
in which
\begin{eqnarray}
	\label{Eq: Matrix def}
	\hat{\textbf{S}}^{(e,\varepsilon)}_{kp} &=& \Big(\frac{d\psi_p }{dx} \,,\, H^{(\varepsilon)}_k(x)   \Big)_{\Omega_e}, \quad e=1,2,\cdots, \varepsilon-1 ,
	\\ \nonumber
	\textbf{S}^{(\varepsilon)}_{kp} &=&
	\,\,
	\Big(\frac{d\psi_p }{dx} \,,\,  \prescript{RL}{x}{\mathcal{D}}_{x_{\varepsilon}}^{\mu} \,
	\Big[  \prescript{(2)}{}{\mathcal{P}}_{k+1}^{\mu}( x ) \, \Big]  \Big)_{\Omega_{\varepsilon}}, 
	\\ \nonumber
	\textbf{M}^{(\varepsilon)}_{kp} &=& \Big( \psi_p ( x  ) \,,\, \,\prescript{(2)}{}{\mathcal{P}}_{k+1}^{\mu}( x ) \Big)_{\Omega_{\varepsilon}},  
	\\ \nonumber
	\textbf{f}^{(\varepsilon)}_k &=& \Big( f \,,\,\prescript{(2)}{}{\mathcal{P}}_{k+1}^{\mu}( x ) \Big)_{\Omega_{\varepsilon}},
\end{eqnarray}
are respectively the \textit{history}, local stiffness, local mass matrices, and local force vector.

%
%%%%%%%%%%%%%%%%%%%%%%%%%%%%%%%%%%%%%%%%%
\subsection{\textbf{Elemental (Local) Operations}: the construction of local matrices $\textbf{S}^{(\varepsilon)}$ and $\textbf{M}^{(\varepsilon)}$, and vector $\textbf{f}^{(\varepsilon)}$}
\label{Sec: lcoal matrices}
%%%%%%%%%%%%%%%%%%%%%%%%%%%%%%%%%%%%%%%%%
%
Here, we provide the analytically obtained expressions of the local stiffness matrix as well as the proper quadrature rules to construct the local mass matrix and force vector in the PG method. 

\begin{figure}[t]
	\center
	\includegraphics[width=0.4\textwidth]{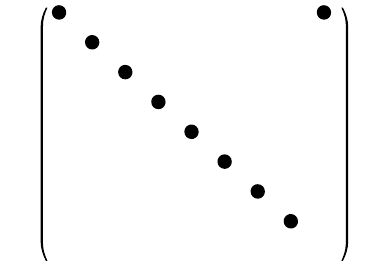}
	\caption{\label{Fig: Stiffness Sparsity} Sparsity of local stiffness matrix}
\end{figure}

\vspace{0.2cm}
\noindent{ \textbf{Elemental (Local) Stiffness Matrix $\textbf{S}^{(\varepsilon)}$} }: given the structure of the basis functions, we first obtain the first ($p=0$) and last column ($p=P$) of the the local stiffness matrix $\textbf{S}^{(\varepsilon)}$, and then, the rest of entries corresponding to the interior modes. Hence,
\begin{eqnarray}
	\label{Eq: Matrix def}
	\textbf{S}^{(\varepsilon)}_{k0} &=& 
	\int_{x_{\varepsilon -1}}^{x_\varepsilon}
	\frac{d\psi_0 }{dx} \, \prescript{RL}{x}{\mathcal{D}}_{x_{\varepsilon}}^{\mu} \,
	\Big[ \prescript{(2)}{}{\mathcal{P}}_{k+1}^{\mu}( x ) \, \Big] 
	dx ,
	\\ \nonumber
	&=&
	\textit{Jac}(\varepsilon , \mu)
	\int_{-1}^{1}
	(\frac{-1 }{2}) (\frac{d\zeta}{dx}) 
	\frac{\Gamma[1+k+\mu]}{\Gamma[1+k]} P_k(\zeta)
	(\frac{dx}{d\zeta})
	d\zeta ,
	\\ \nonumber
	&=&
	-\textit{Jac}(\varepsilon , \mu)
	\frac{\Gamma(1+k+\mu)}{2 \,\, \Gamma(1+k)}
	\int_{-1}^{1}
	P_k(\zeta)
	d\zeta ,
	\\ \nonumber
	&=&
	-
	\textit{Jac}(\varepsilon , \mu) 
	\frac{\Gamma(1+k+\mu)}{\Gamma(1+k)}
	\delta_{k , 0}, \quad (\textit{by the orthogonality})
\end{eqnarray}
in which the Jacobian constant, associated with the element $\varepsilon$ and the fractional order $\mu$, is $\textit{Jac}(\varepsilon , \mu) = (\frac{2}{x_{\varepsilon} - x_{\varepsilon-1} })^{\mu}$; hence the first column of the local stiffness matrix for $\varepsilon = 1,2,\cdots, N_{el}$ is obtained as
\begin{equation}
	\boxed{
		\textbf{S}^{(\varepsilon)}_{k0} =
		-
		(\frac{2}{x_{\varepsilon} - x_{\varepsilon-1} })^{\mu}
		\frac{\Gamma(1+k+\mu)}{\Gamma(1+k)}
		\delta_{k , 0}, \quad k=0,1,\cdots, P.
	}
\end{equation}
Similarly, we can obtain the last column of the local stiffness matrix $\textbf{S}^{(\varepsilon)}_{kP}$ as
\begin{equation}
	\boxed{
		\textbf{S}^{(\varepsilon)}_{kP} =
		(\frac{2}{x_{\varepsilon} - x_{\varepsilon-1} })^{\mu}
		\frac{\Gamma(1+k+\mu)}{\Gamma(1+k)}
		\delta_{k , 0} = - \textbf{S}^{(\varepsilon)}_{k0}, \quad k=0,1,\cdots, P.
	}
\end{equation}
In order to obtain the rest of entries of $\textbf{S}^{(\varepsilon)}_{kp}$ ($k=0,1, \cdots, P$ and $p= 1,2,\cdots, P-1$), we carry out the integration-by-parts and transfer another derivative onto the test function, taking into account that the interior modes vanish at the boundary points $x_{\varepsilon}$ and $x_{\varepsilon -1}$. Therefore,

\begin{eqnarray}
	\label{Eq: Matrix calculation stiffness}
	\textbf{S}^{(\varepsilon)}_{kp} &=& 
	\int_{x_{\varepsilon -1}}^{x_\varepsilon}
	\frac{d \psi_p}{dx} \, \prescript{RL}{x}{\mathcal{D}}_{x_{\varepsilon}}^{\mu} \,
	\Big[  \prescript{(2)}{}{\mathcal{P}}_{k+1}^{\mu}( x ) \, \Big] 
	dx ,
	\\ \nonumber
	&=&
	-\int_{x_{\varepsilon -1}}^{x_\varepsilon}
	\psi_p (x)\,\,\,\frac{d}{dx} \prescript{RL}{x}{\mathcal{D}}_{x_{\varepsilon}}^{\mu} \,
	\Big[ \prescript{(2)}{}{\mathcal{P}}_{k+1}^{\mu}( x ) \, \Big] 
	dx ,
	\\ \nonumber
	&=&
	-\int_{-1}^{1}
	\psi_p (\zeta)\,\,\,\frac{d}{d\zeta} \frac{d\zeta}{dx} \,\,
	%
	%  (\frac{2}{x_{\varepsilon} - x_{\varepsilon -1}})^{\mu}
	\textit{Jac}(\varepsilon , \mu)
	\prescript{RL}{\zeta}{\mathcal{D}}_{1}^{\mu} \,
	\Big[  \prescript{(2)}{}{\mathcal{P}}_{k+1}^{\mu}( \zeta  ) \, \Big] 
	\frac{dx}{d\zeta}
	d\zeta ,
	\\ \nonumber
	&=&
	%
	% - (\frac{2}{x_{\varepsilon} - x_{\varepsilon -1}})^{\mu}
	- \textit{Jac}(\varepsilon , \mu)
	\frac{\Gamma(1+k+\mu)}{4 \,\, \Gamma(1+k)}
	\int_{-1}^{1}
	(1- \zeta)(1+ \zeta) P^{1,1}_{p-1}(\zeta) \,\,
	\,
	\frac{d}{d\zeta} 
	\Big[ 
	P_k(\zeta) \Big] 
	d\zeta , 
	\\ \nonumber
	&=&
	%
	% - (\frac{2}{x_{\varepsilon} - x_{\varepsilon -1}})^{\mu}
	- \textit{Jac}(\varepsilon , \mu)
	\frac{\Gamma(1+k+\mu)}{4 \,\, \Gamma(1+k)}
	\frac{k+1}{2}
	\int_{-1}^{1}
	(1- \zeta)(1+ \zeta) P^{1,1}_{p-1}(\zeta) \,\,
	\,
	P^{1,1}_{k-1}(\zeta) 
	d\zeta.
\end{eqnarray}
Hence, for $\varepsilon = 1,2,\cdots, N_{el}$,

\begin{equation}
	\boxed{
		\textbf{S}^{(\varepsilon)}_{kp} =
		-(\frac{2}{x_{\varepsilon} - x_{\varepsilon -1}})^{\mu}
		\frac{\Gamma(1+k+\mu) (k+1)}{8 \,\, \Gamma(1+k)}
		C^{1,1}_{k-1} \delta_{k,p},
		\quad
		\begin{aligned}
			& k=0, \cdots, P, 
			\\
			&p= 1,\cdots, P-1,
		\end{aligned}
	}
\end{equation}
where $C^{1,1}_{k-1} $ represents the corresponding orthogonality constant of Jacobi polynomials of order $k-1$ with parameters $\alpha = \beta =1$. We note that the entries of $\textbf{S}^{(\varepsilon)}_{kp}$ are obtained analytically using the orthogonality of Jacobi polynomial. Also, the interior modes lead to a \textit{diagonal} matrix due to $\delta_{k,p}$. Fig. \ref{Fig: Stiffness Sparsity} shows the sparsity of the local stiffness matrix.

%The rest of entries of $\textbf{S}^{(\varepsilon)}_{kp}$ for $k=0,1, \cdots, P$ and $p= 1,2,\cdots, P-1$ are also obtained analytically, since the interior modes lead to a \textit{diagonal} matrix. It is done by carrying out the integration-by-parts once more to transfer another derivative onto the test function, taking into account that the bubble modes vanish at the boundary points $x_{\varepsilon}$ and $x_{\varepsilon -1}$. So,

%
%******************************************************************************************

\vspace{0.2cm}
\noindent{ \textbf{Elemental (Local) Mass Matrix $\textbf{M}^{(\varepsilon)}$} }: similarly, we first obtain the corresponding first ($p=0$) and last column ($p=P$) of the local mass matrix $\textbf{M}^{(\varepsilon)}$, and then, we compute the rest of entries associated with the interior modes.
\begin{eqnarray}
	\nonumber
	\textbf{M}^{(\varepsilon)}_{k0} &=& 
	\int_{x_{\varepsilon-1}}^{x^{\varepsilon }}
	\psi_0 (x) \, \,\prescript{(2)}{}{\mathcal{P}}_{k+1}^{\mu}(x) dx ,
	\\ \nonumber
	&=&  
	(\frac{x_{\varepsilon}  -  x_{\varepsilon -1}}{2})
	\int_{-1}^{1}
	\psi_0 (\zeta) \, \,\prescript{(2)}{}{\mathcal{P}}_{k+1}^{\mu}( \zeta ) d\zeta ,
	\\ \nonumber
	&=&  
	\frac{1}{2}
	(\frac{x_{\varepsilon}  -  x_{\varepsilon -1}}{2})
	\int_{-1}^{1}
	(1 - \zeta)^{1+\mu} \, P^{\mu , -\mu}_{k} ( \zeta ) d\zeta.
\end{eqnarray}
Hence, for $\varepsilon = 1,2,\cdots, N_{el}$,
\begin{equation}
	\nonumber
	\boxed{
		\textbf{M}^{(\varepsilon)}_{k0}  =
		(\frac{x_{\varepsilon}  -  x_{\varepsilon -1}}{4})
		\sum_{q=1}^{Q} w^{ 1+\mu , 0}_q  \,\, P^{\mu , -\mu}_{k} (\, z^{1+\mu , 0}_q \,)
		, \quad k=0,1,\cdots, P ,
	}
\end{equation}
where $\{ w^{ 1+\mu , 0}_q  \,,\, z^{1+\mu , 0}_q  \}_{q=1}^{Q}$ are the Gauss-Lobatto-Jacobi weights and points corresponding to the parameters $\alpha = 1+\mu$ and $\beta =0$. By similar steps, we obtain 
\begin{eqnarray}
	\nonumber
	\textbf{M}^{(\varepsilon)}_{kP} &=& 
	\int_{x_{\varepsilon -1}}^{x^{\varepsilon }}
	\psi_P (x) \, \,\prescript{(2)}{}{\mathcal{P}}_{k+1}^{\mu}(x) dx ,
	\\ \nonumber
	&=&  
	(\frac{x_{\varepsilon}  -  x_{\varepsilon -1}}{2})
	\int_{-1}^{1}
	\psi_P ( \zeta  ) \, \,\prescript{(2)}{}{\mathcal{P}}_{k+1}^{\mu}( \zeta ) d\zeta ,
	\\ \nonumber
	&=&  
	\frac{1}{2}
	(\frac{x_{\varepsilon}  -  x_{\varepsilon -1}}{2})
	\int_{-1}^{1}
	(1 + \zeta) 
	(1 - \zeta)^{\mu} \, P^{\mu , -\mu}_{k} ( \zeta ) d\zeta.
\end{eqnarray}
Therefore, for $\varepsilon = 1,2,\cdots, N_{el}$,
\begin{equation}
	\nonumber
	\boxed{
		\textbf{M}^{(\varepsilon)}_{kP}  =
		(\frac{x_{\varepsilon}  -  x_{\varepsilon -1}}{4})
		\sum_{q=1}^{Q} w^{ \mu , 1}_q  \,\, P^{\mu , -\mu}_{k} (\, z^{\mu , 1}_q \,) , \quad k=0,1,\cdots, P ,
	}
\end{equation}
where $\{ w^{\mu , 1}_q  \,,\, z^{\mu , 1}_q  \}_{q=1}^{Q}$ are the Gauss-Lobatto-Jacobi weights and points corresponding to the parameters $\alpha = \mu$ and $\beta = 1$. The rest of the entries of the local mass matrix are then obtained as
\begin{eqnarray}
	\nonumber
	\textbf{M}^{(\varepsilon)}_{kp} &=& 
	\int_{x_{\varepsilon -1}}^{x^{\varepsilon }}
	\psi_p (x) \, \,\prescript{(2)}{}{\mathcal{P}}_{k+1}^{\mu}(x) dx ,
	\\ \nonumber
	&=&  
	(\frac{x_{\varepsilon}  -  x_{\varepsilon -1}}{2})
	\int_{-1}^{1}
	\psi_p ( \zeta  ) \, \,\prescript{(2)}{}{\mathcal{P}}_{k+1}^{\mu}( \zeta ) d\zeta ,
	\\ \nonumber
	&=&  
	\frac{1}{4}
	(\frac{x_{\varepsilon}  -  x_{\varepsilon -1}}{2})
	\int_{-1}^{1}
	(1 + \zeta) 
	(1 - \zeta)^{1+\mu} \, P^{\mu , -\mu}_{k} ( \zeta ) d\zeta, 
\end{eqnarray}
and thus, for $\varepsilon = 1,2,\cdots, N_{el}$,
\begin{equation}
	\nonumber
	\boxed{
		\textbf{M}^{(\varepsilon)}_{kp}  =
		(\frac{x_{\varepsilon}  -  x_{\varepsilon -1}}{8})
		\sum_{q=1}^{Q} w^{ 1+\mu , 1}_q  \,\, P^{\mu , -\mu}_{k} (\, z^{1+\mu , 1}_q \,)  , \quad k=0,1,\cdots, P ,
	}
\end{equation}
where $\{ w^{1+\mu , 1}_q  \,,\, z^{1+\mu , 1}_q  \}_{q=1}^{Q}$ are the Gauss-Lobatto-Jacobi weights and points corresponding to the parameters $\alpha =1+ \mu$ and $\beta = 1$.

\vspace{0.2cm}
\noindent{ \textbf{Elemental (Local) Load Vector $\textbf{f}^{(\varepsilon)}$} }: the local load vector is obtained as:

\begin{eqnarray}
	\nonumber
	\textbf{f}^{(\varepsilon)}_k =
	\int_{x_{\varepsilon -1}}^{x^{\varepsilon }}
	f(x) \,\, \prescript{(2)}{}{\mathcal{P}}_{k+1}^{\mu}( x  ) 
	dx 
	=
	(\frac{x_{\varepsilon}  -  x_{\varepsilon -1}}{2})
	\int_{-1}^{1}
	(1-\zeta)^{\mu} \,\, 
	f(x(\zeta)) \,\, P^{\mu , -\mu}_{k}(\zeta) \,\,  
	%\frac{dx^{\varepsilon}}{d\zeta} 
	d\zeta.
\end{eqnarray}
Hence, for $\varepsilon = 1,2,\cdots, N_{el}$,
\begin{equation}
	\nonumber
	\boxed{
		\textbf{f}^{(\varepsilon)}_k  =
		(\frac{x_{\varepsilon}  -  x_{\varepsilon -1}}{2})
		\sum_{q=1}^{Q} w^{ \mu , 0}_q  \,\,f(\,  x^{\varepsilon}(\zeta_q)  \,)\, P^{\mu , -\mu}_{k} (\, z^{\mu , 0}_q \,) , \quad k=0,1,\cdots, P ,
	}
\end{equation}
where $\{ w^{\mu , 0}_q  \,,\, z^{\mu , 0}_q  \}_{q=1}^{Q}$ are the Gauss-Lobatto-Jacobi weights and points corresponding to the parameters $\alpha = \mu$ and $\beta = 0$.

%
%%%%%%%%%%%%%%%%%%%%%%%%%%%%%%%%%%%%%%%%%
\subsection{ \textbf{Non-Local Operation}: the construction of history matrix  $\hat{\textbf{S}}^{(e,\varepsilon)}$}
\label{Sec: }
%%%%%%%%%%%%%%%%%%%%%%%%%%%%%%%%%%%%%%%%%
%
The most challenging part of constructing the linear system is to compute the global history matrix $\hat{\textbf{S}}^{(e,\varepsilon)}$. The history matrix relates the current element $\varepsilon = 1,2,\cdots,N_{el}$ to its past elements $e = 1,2,\cdots,\varepsilon-1$ by
\begin{eqnarray}
	\label{Eq: global history mat}
	\hat{\textbf{S}}^{(e,\varepsilon)}_{kp} &=& 
	\int_{x_{e-1}}^{x_{e}}
	\frac{d\psi_p }{dx} \, H^{(\varepsilon)}_k(x) \,\,  dx,
	\quad k=0, \cdots, P, \quad p= 1,\cdots, P-1,
\end{eqnarray}
where $H^{(\varepsilon)}_k(x)$ is given in \eqref{Eq: frac der test function-2} as
\begin{align*}
	H_k^{(\varepsilon)}(x) 
	& = \frac{-1}{\Gamma(1-\mu)} \frac{d}{dx} \int_{x_{\varepsilon -1}}^{x^{\varepsilon}}  \frac{\prescript{(2)}{}{\mathcal{P}}_{k+1}^{\mu}(s)  }{(s - x)^{\mu}} ds ,
	\\
	& = \frac{-\mu}{\Gamma(1-\mu)} \int_{x_{\varepsilon -1}}^{x^{\varepsilon}}  \frac{\prescript{(2)}{}{\mathcal{P}}_{k+1}^{\mu}(s)  }{(s - x)^{1+\mu}} ds ,
\end{align*}
in which, $x \in \Omega_e = \left[ x_{e-1} , x_{e} \right]$ and $s \in \Omega_\varepsilon = \left[ x_{\varepsilon-1} , x_{\varepsilon} \right]$. By performing the following affine mappings
\begin{eqnarray}
	\nonumber
	s &=& \frac{x_{\varepsilon}   +  x_{\varepsilon-1}  }{2}  + \frac{x_{\varepsilon}   -  x_{\varepsilon-1}  }{2}  \zeta ,
	\\ \nonumber
	x &=& \frac{x_{e}   +  x_{e-1}  }{2}  + \frac{x_{e}   -  x_{e-1}  }{2}  \xi ,
\end{eqnarray}
from $\Omega_e$ and $\Omega_{\varepsilon}$ to the standard element $\left[ -1 , 1 \right]$, the history function $H_k^{(\varepsilon , e)}(\xi) = H_k^{(\varepsilon)}(x)$ is obtained as 
\begin{align}
	\label{Eq: History function general case}
	&H_k^{(\varepsilon , e) }(\xi)
	%&=&
	%%
	%(  \frac{x_{\varepsilon}   -  x_{\varepsilon-1}  }{2}   )
	%\frac{-\mu}{\Gamma(1-\mu)} \int_{-1}^{1}  \frac{\prescript{(2)}{}{\mathcal{P}}_{k+1}^{\mu}(\zeta  )  }{[s^{\varepsilon}(\zeta)- x(\xi)]^{1+\mu}}  d\zeta,
	%%
	\\ \nonumber
	&=
	(  \frac{x_{\varepsilon}   -  x_{\varepsilon-1}  }{2}   )
	\frac{-\mu}{\Gamma(1-\mu)} \int_{-1}^{1}  \frac{\prescript{(2)}{}{\mathcal{P}}_{k+1}^{\mu}(\zeta  )    d\zeta}{
		\Big[      \frac{(x_{\varepsilon}   +  x_{\varepsilon-1}) - (x_{e}   +  x_{e-1})  }{2}   + \frac{x_{\varepsilon}   -  x_{\varepsilon-1}  }{2}  \zeta - \frac{x_{e}   -  x_{e-1}  }{2}  \xi          \Big]^{1+\mu}}.
\end{align}
%
%and $\forall x \in \Omega_e,\,   e \leq \varepsilon-1$,
%%
%\begin{eqnarray}
%%\label{Eq: history-func}
%\nonumber
%%
%H_k^{(\varepsilon)}(x) &=&
%%
%\frac{-1}{\Gamma(1-\mu)} \frac{d}{dx} \int_{0}^{x^{\varepsilon}}  \frac{v^{\varepsilon}(x) }{(s - x)^{\mu}} ds = \frac{-1}{\Gamma(1-\mu)} \sum_{e=1}^{\varepsilon}  \frac{d}{dx} \int_{x_{e-1}}^{x_{e}}  \frac{v^{\varepsilon}(x) }{(s - x)^{\mu}} ds
%%
%\\ \nonumber
%%
%&=&
%%
%\frac{-1}{\Gamma(1-\mu)} \sum_{e=1}^{\varepsilon -1}  \frac{d}{dx} \int_{x_{e-1}}^{x_{e}}  \frac{ 0 }{(s - x)^{\mu}} ds
%+
% \frac{-1}{\Gamma(1-\mu)} \frac{d}{dx} \int_{x_{\varepsilon -1}}^{x^{\varepsilon}}  \frac{\prescript{(2)}{}{\mathcal{P}}_{k+1}^{\mu}( s^e(\zeta) )  }{(s- x)^{\mu}} ds
% %
%\\ \nonumber
%%
%&=&
%%
% \frac{-\mu}{\Gamma(1-\mu)} \int_{x_{\varepsilon -1}}^{x^{\varepsilon}}  \frac{\prescript{(2)}{}{\mathcal{P}}_{k+1}^{\mu}( s^e(\zeta)  )  }{(s- x)^{1+\mu}} ds.
%%
%\end{eqnarray}
%%
If the mesh is ``uniform'', then 
\begin{align}
	\nonumber
	&
	\frac{x_{\varepsilon}   +  x_{\varepsilon-1}  }{2}  = \frac{2\varepsilon -1}{2} \Delta x , \quad
	\frac{x_{e}   +  x_{e-1}  }{2}  = \frac{2e -1}{2} \Delta x ,
	\\
	\label{Eq: uniform mesh}
	&
	\frac{ x_{\varepsilon} -  x_{\varepsilon -1} }{2} = \frac{ x_{e} -  x_{e -1} }{2} = \frac{\Delta x}{2} ,
\end{align}
and thus,
\begin{eqnarray}
	\nonumber
	H_k^{(\varepsilon , e) }(\xi) &=&
	\frac{\Delta x  }{2}   
	\frac{-\mu}{\Gamma(1-\mu)} \int_{-1}^{1} 
	\frac{\prescript{(2)}{}{\mathcal{P}}_{k+1}^{\mu}(\zeta  )  }{
		\Big[      \frac{(2\varepsilon -1) - (2e -1)  }{2}\Delta x     + \frac{\Delta x  }{2}  (\zeta -  \xi)          \Big]^{1+\mu}} d\zeta,
	\\ \nonumber
	&=&
	\frac{-\mu}{\Gamma(1-\mu)} 
	(\frac{2}{\Delta x})^{\mu}
	\,\,
	\int_{-1}^{1} 
	\frac{\prescript{(2)}{}{\mathcal{P}}_{k+1}^{\mu}(\zeta  )  }{
		\Big[ 2(\varepsilon - e)  +  \zeta -  \xi \Big]^{1+\mu}} d\zeta,
	\\ \label{Eq: History function unifrom mesh}
	&=&
	\frac{-\mu}{\Gamma(1-\mu)} 
	(\frac{2}{\Delta x})^{\mu}
	\,\,
	\int_{-1}^{1} 
	\frac{\prescript{(2)}{}{\mathcal{P}}_{k+1}^{\mu}(\zeta  )  }{
		\Big[ 2 \, \Delta \varepsilon  +  \zeta -  \xi \Big]^{1+\mu}} d\zeta,
\end{eqnarray}
%
%$\forall \xi = \xi^{e} \in [-1,1], k=0,1,\cdots, P$.
%%
%\begin{eqnarray}
%\label{Eq: history-func h uniform}
%h^{(\varepsilon , e)}_k(\xi) \equiv \int_{-1}^{1} 
%%
% \frac{\prescript{(2)}{}{\mathcal{P}}_{k+1}^{\mu}(\zeta  )  }{
%\Big[ 2(\varepsilon - e)  +  \zeta -  \xi \Big]^{1+\mu}} d\zeta, \quad \forall \xi = \xi^{e} \in [-1,1], k=0,1,\cdots, P.
%%
%\end{eqnarray}
%%
where $\Delta \varepsilon = \varepsilon - e > 0$, denotes the element difference between the current element $\varepsilon$ and the $e$-th element. Next, we expand the poly-fractonomials $\prescript{(2)}{}{\mathcal{P}}_{k+1}^{\mu}(\zeta  )$ in terms of fractonomials $(1-\zeta)^{\mu + m}$ as
\begin{eqnarray}
	\label{Eq: poly-fractonomial expansion}
	\prescript{(2)}{}{\mathcal{P}}_{k+1}^{\mu}(\zeta  ) = (1-\zeta)^{\mu} P^{\mu , -\mu}_{k}(\zeta)
	=
	\sum_{m=0}^{k} C_{km} (1-\zeta)^{\mu + m},
\end{eqnarray}
in which $C_{km} = {k+m \choose m} {k+\mu \choose k- m} (-\frac{1}{2})^m$ is a lower-triangle matrix. Therefore, \eqref{Eq: History function unifrom mesh} can be written as
\begin{eqnarray}
	\label{Eq: history-func H uniform-2}
	H^{(\varepsilon , e)}_k(\xi) = 
	\frac{-\mu}{\Gamma(1-\mu)} 
	(\frac{2}{\Delta x})^{\mu}
	\sum_{m=0}^{k} C_{km} 
	h^{(\varepsilon , e)}_m(\xi),
\end{eqnarray}
where we call
\begin{align}
	\label{Eq: modal memory mode}
	h^{(\varepsilon , e)}_m(\xi) \equiv \int_{-1}^{1} 
	\frac{(1-\zeta)^{\mu + m}  }{[ 2 \, \Delta \varepsilon  +  \zeta -  \xi ]^{1+\mu}} d\zeta, \quad m=0,1, \cdots, k,
\end{align}
the \textit{(modal) memory mode}. Also, $h^{(\varepsilon , e)}_m(\xi)$ can be obtained analytically as
\begin{eqnarray}
	\label{Eq: history-func h uniform}
	h^{(\varepsilon , e)}_m(\xi)  = \frac{2( \Delta \varepsilon - \xi/2)}{1+m +\mu} 
	\Big[ h_{m, I}(\xi , \Delta \varepsilon)  + h_{m, II}(\xi , \Delta \varepsilon)  + h_{m, III}(\xi , \Delta \varepsilon) 
	\Big],
\end{eqnarray}
in which
%%
%\begin{eqnarray}
%\label{Eq: h-I uniform}
%%
% h_{m, I}(\xi , \Delta \varepsilon) &=& - \frac{ _2F_1\Big(  1 \,,\, 1+m \,,\, 2+m+\mu \,,\,  \frac{1}{1+2\Delta \varepsilon - \xi} \Big)     }{    1 +2\Delta \varepsilon - \xi   }
%%
%\\ \nonumber
%%
%  h_{m, II}(\xi , \Delta \varepsilon) &=& 
%  %
%  -    \Big(   \frac{1-2 \Delta \varepsilon  + \xi   }{ -2\Delta \varepsilon  + \xi  }      \Big)^{-\mu}\,  (2^{1+m+\mu} )  \,   
%  %
%   \frac{ _2F_1\Big(  1 \,,\, 1+m \,,\, 2+m+\mu \,,\,  \frac{2}{-1 -2\Delta \varepsilon + \xi} \Big)     }{    -1 - 2\Delta \varepsilon + \xi   },
%   %
%\\ \nonumber
%%
%  h_{m, III}(\xi , \Delta \varepsilon) &=&   
%  %
%   \frac{_2F_1\Big(  1 \,,\, 1+m \,,\, 2+m+\mu \,,\,  \frac{1}{ -2\Delta \varepsilon + \xi} \Big)     }{    - 2\Delta \varepsilon + \xi   }
%%
%\end{eqnarray}
%%
%
\begin{align}
	\label{Eq: h-I uniform}
	h_{m, I}(\xi , \Delta \varepsilon) &= - Z_I(\xi , \Delta \varepsilon ) \,\, _2F_1\Big(  1 \,,\, 1+m \,,\, 2+m+\mu \,,\,  Z_I(\xi, \Delta \varepsilon)  \Big)   , 
	\\ \nonumber
	h_{m, II}(\xi , \Delta \varepsilon) &= 
	\Big(   \frac{1-2 \Delta \varepsilon  + \xi   }{ -2\Delta \varepsilon  + \xi  }      \Big)^{-\mu}  (2^{m+\mu} )  \,\,   
	Z_{II}(\xi, \Delta \varepsilon) 
	_2F_1\Big(  1 \,,\, 1+m \,,\, 2+m+\mu \,,\,  Z_{II}(\xi, \Delta \varepsilon)  \Big)    ,
	\\ \nonumber
	h_{m, III}(\xi , \Delta \varepsilon) &=   
	- Z_{III}(\xi, \Delta \varepsilon) _2F_1\Big(  1 \,,\, 1+m \,,\, 2+m+\mu \,,\,  Z_{III}(\xi, \Delta \varepsilon) \Big)    ,
\end{align}
and the group variables are $Z_I(\xi, \Delta \varepsilon)  = \frac{1}{1 +2\Delta \varepsilon - \xi} $, $Z_{II}(\xi, \Delta \varepsilon) = \frac{-2}{-1 - 2\Delta \varepsilon + \xi }  $, and $Z_{III}(\xi, \Delta \varepsilon)=\frac{1 }{    - 2\Delta \varepsilon + \xi   }$. 

Moreover, The derivative of the basis function in the standard element is given by
\begin{eqnarray}
	\label{Eq: basis derivative}
	\frac{d\psi_p (\zeta)}{d\zeta} =
	\begin{cases}
		\frac{-1}{2},\quad & p=0, \\
		%\\
		\frac{d}{d\zeta} \Big[  (\frac{1- \zeta}{2})(\frac{1+ \zeta}{2}) P^{1,1}_{p-1}(\zeta) \Big],\quad & p=1,2,\cdots, P-1, \\
		%\\
		\frac{1}{2},\quad & p=P.
	\end{cases}
\end{eqnarray}
Therefore, by \eqref{Eq: poly-fractonomial expansion} and \eqref{Eq: basis derivative}, the entries of the history matrix can be efficiently computed using a Gauss quadrature. Hence:
\begin{equation}
	\label{Eq: global history mat-standard domain}
	\boxed{
		\hat{\textbf{S}}^{(\varepsilon , e)}_{kp} \equiv
		\hat{\textbf{S}}^{(\Delta \varepsilon)}_{kp} =
		\int_{-1}^{1}
		\frac{d\psi_p }{d\xi} \, H_k(\xi , \Delta \varepsilon)   d\xi, \quad k,p =0,1,\cdots, P.
	}
\end{equation}

\vspace{0.25 cm}
\begin{rem}
	\label{Rem: on the uniform mesh and history matrix}
	We note that when a uniform mesh is employed, the history function $H^{(\varepsilon , e)}_{k} (\xi) \equiv H_{k}(\xi , \Delta \varepsilon)$, defined in the standard element, only depends on the ``element difference'', $\Delta \varepsilon = \varepsilon - e$. This is significant since one only needs to construct $N_{el}-1$ history function, and thus, history matrices $\hat{\textbf{S}  }^{(e , \varepsilon)}$. 
\end{rem}
\vspace{0.25 cm}
%

%
%%%%%%%%%%%%%%%%%%%%%%%%%%%%%%%%%%%%%%%%%
\subsection{\textbf{Assembling the Global System with Local Test Functions}}
\label{Sec: Assembling Local Test}
%%%%%%%%%%%%%%%%%%%%%%%%%%%%%%%%%%%%%%%%%
%
We generalize the notion of global linear system assembly by taking into account the presence of the history stiffness matrices and recalling that the corresponding local mass matrix $\textbf{M}^{(\varepsilon)}$ or the local load-vector $\textbf{f}^{(\varepsilon)}$ do not contribute to any history calculations. We impose the $C^0-continuity$ by employing the ``mapping arrays'', map[e][p], defined as
\begin{equation}
	\label{Eq: mapping array}
	map[e][p] = P(e-1) + p, \quad p=1,2,\cdots, P, \quad e= 1, 2, \cdots, N_{el},
\end{equation}
as for instance in Mathematica, the first entry of a vector is labelled by 1 rather than 0 as in C++. Then, the corresponding $(P+1)\times(P+1)$ ``local" linear system, which is associated with the element $\Omega_{\varepsilon}$, is obtained as
\begin{equation}
	\label{Eq: local linear system}
	\mathcal{M}^{(\varepsilon)} = \textbf{S}^{(\varepsilon)} - \lambda \textbf{M}^{(\varepsilon)}.
\end{equation}
We assemble the corresponding global linear matrix $\mathbb{M}_G$ and the global load-vector $\mathbb{F}_G$ as follows:
\begin{mymathbox}
	\vspace{-0.25 in}
	\begin{eqnarray}
		\nonumber
		\textit{do}\,\,\varepsilon =1, N_{el} 
		\hspace{11cm} 
		\\ \nonumber
		\textit{do}\,\,k=1, P+1
		\hspace{10cm} 
		\\ \nonumber
		\mathbb{F}_G\Big[ \,\, map[\varepsilon][k]  \,\,   \Big] = \textbf{f}^{(\varepsilon)} [k]
		\hspace{8cm} 
		\\ \nonumber
		\textit{do}\,\,p=1, P+1
		\hspace{9.5cm} 
		\\ \nonumber
		\mathbb{M}_G \Big[\,  map[\varepsilon][k]      \,\Big]\Big[\,   map[\varepsilon][p]  \,\Big]  = \mathbb{M}_G \Big[\,  map[\varepsilon][k]      \,\Big]\Big[\,   map[\varepsilon][p]  \,\Big] + \mathcal{M}^{(\varepsilon)}[k][p]
		\hspace{-1cm} 
		\\ \nonumber
		{\color{blue}
			\textit{do}\,\,e=1, \varepsilon-1
			\hspace{9cm} 
		}
		\\ \nonumber
		{\color{blue}
			\mathbb{M}_G \Big[\,  map[  {\color{red}\varepsilon}  ][k]      \,\Big]\Big[\,   map[  {\color{red}e}  ][p]  \,\Big]  = \mathbb{M}_G \Big[\,  map[  {\color{red}\varepsilon}  ][k]      \,\Big]\Big[\,   map[  {\color{red}e}  ][p]  \,\Big]  + \hat{\textbf{S}}^{(\Delta \varepsilon)}[k][p]
			\hspace{-1.3cm} 
		}
		%\\ \nonumber
		%{\color{blue}
		%End
		%}
		%\hspace{10.60 cm} 
		%\\ \nonumber
		%End
		%\hspace{11.2 cm} 
		%\\ \nonumber
		%End
		%\hspace{11.60 cm} 
		\\ \nonumber
		End
		\hspace{12.15 cm} 
	\end{eqnarray}
\end{mymathbox}
This global operation leads to the following linear system:
\begin{equation}
	\label{Eq: global linear system}
	\boxed{
		\mathbb{M}_G \,\, \hat{u}_G = \mathbb{F}_G,
	}
\end{equation}
in which $\hat{u}_G$ denotes the global degrees of freedom. The homogeneous Dirichlet boundary conditions are enforced by ignoring the first and the last rows also the first and the last columns of the global matrix, in addition to ignoring the first and last entries of the load matrix. We do so since we already know that $\hat{u}^{1}_0  = \hat{u}^{N_{el}}_P = 0$. 
%
%******************************************************************************************
\begin{figure}[t]
	\center
	\includegraphics[width=0.65\textwidth]{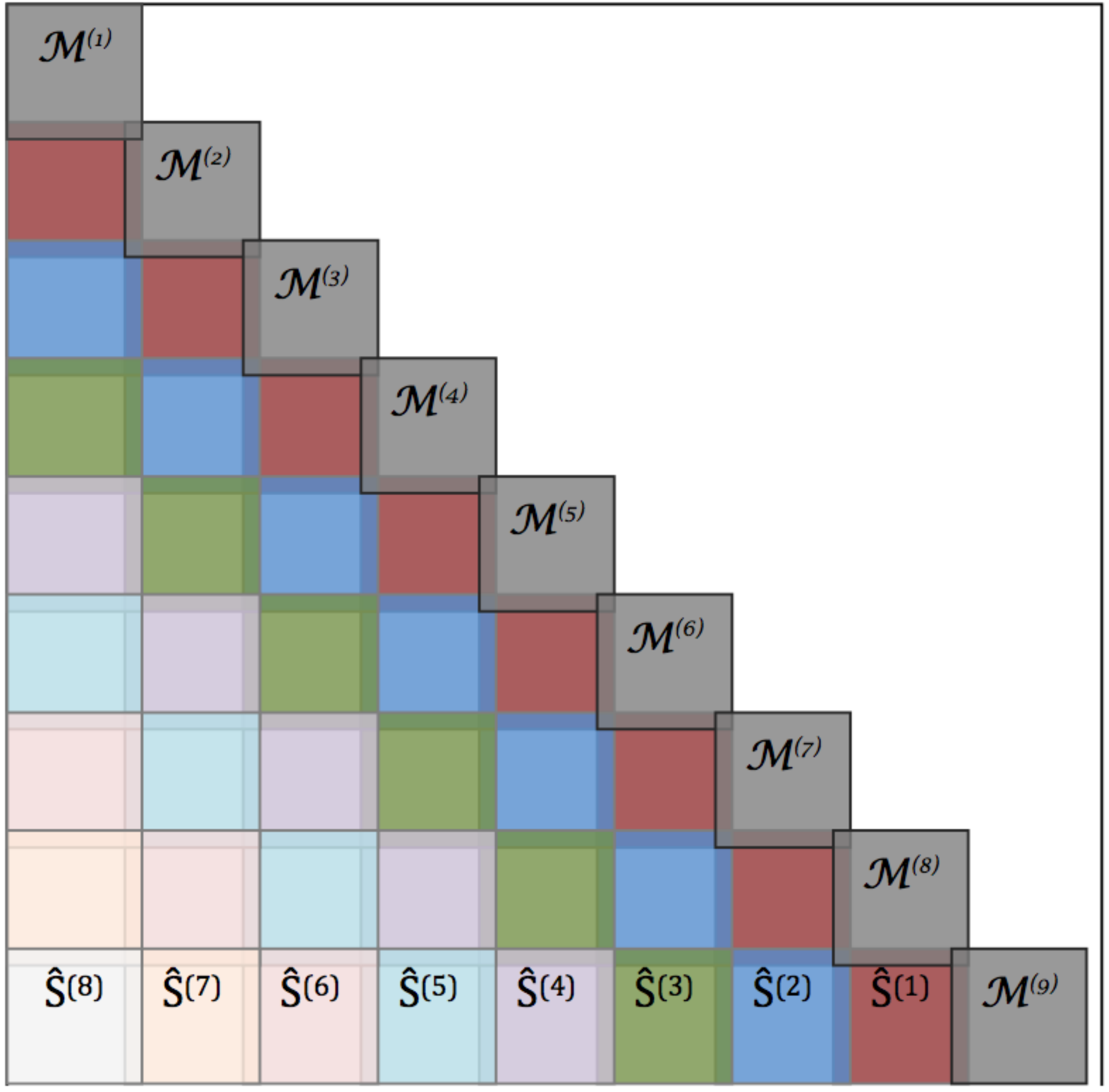}
	\caption{\label{Fig: Global Mat}
		The assembled global matrix corresponding to a uniform grid with $N_{el} = 9$. In this global matrix, $\mathcal{M}^{(\varepsilon)} = \textbf{S}^{(\varepsilon)} - \lambda \textbf{M}^{(\varepsilon)}$, $\varepsilon =1,2,\cdots, N_{el}$, represents the local matrix, associated with the element $\Omega_{\varepsilon}$. To fill the lower-triangular block matrices, we construct only $(N_{el}-1)$ history matrices $\hat{ \textbf{S}   }^{( \Delta \varepsilon )}$, where $\Delta \varepsilon =1,2, .., N_{el} -1$, rather than $\frac{N_{el}(N_{el} -1)}{2}$ matrices.
	}
\end{figure}
% 

%
%%%%%%%%%%%%%%%%%%%%%%%%%%%%%%%%%%%%%%%%%
\subsection{\textbf{Scattering from the Global to Local Degrees of Freedom}}
\label{Sec: Scatterign}
%%%%%%%%%%%%%%%%%%%%%%%%%%%%%%%%%%%%%%%%%
%
Once again, due to our $C^0$-continuity and the decomposition of our basis functions into boundary and interior modes, we have $\hat{u}^{e-1}_P = \hat{u}^{e}_0$.  That leads to the following standard scattering process from the global to local degrees of freedom (see e.g., \cite{Karniadakis2005}):
\begin{mymathbox2}
	\vspace{-0.25 in}
	\begin{eqnarray}
		\nonumber
		\textit{do}\,\,\varepsilon =1, N_{el} 
		\hspace{11cm} 
		\\ \nonumber
		\textit{do}\,\,k=1, P+1
		\hspace{10cm} 
		\\ \nonumber
		\hat{u}^{\varepsilon} [k] = \hat{u}_G [\,\, map[\varepsilon][k]     \,\,]
		\hspace{8 cm} 
		%\\ \nonumber
		%End
		%\hspace{11.60 cm} 
		\\ \nonumber
		End
		\hspace{12.4 cm} 
	\end{eqnarray}
\end{mymathbox2}
%

%\newpage

%
%%%%%%%%%%%%%%%%%%%%%%%%%%%%%%%%%%%%%%%%%
\subsection{\textbf{Off-Line Computation of History Matrices and History Retrieval}}
\label{Sec: offline computation}
%%%%%%%%%%%%%%%%%%%%%%%%%%%%%%%%%%%%%%%%%
%

As mentioned in remark \ref{Rem: on the uniform mesh and history matrix} (on uniform grid generation), the history matrices solely depend on the element difference, $\Delta \varepsilon = \varepsilon - e$. Thus, for all local elements $\varepsilon$, where $\varepsilon =1,2,\cdots, N_{el} $, the history matrices corresponding to the past element $e$ with similar element difference, are the same. See Fig. \ref{Fig: Global Mat}, where similarly-colored blocks represent the same history matrix and one can see that, for example, all the history matrices adjacent to the local stiffness matrices have the same element difference, $\Delta \varepsilon = 1$, and thus are in the same color. Therefore, given number of element $N_{el}$, we only need to construct the total number of $N_{el} - 1$ history matrices.

For a maximum number of elements, $N_{el}|_{max}$, and a maximum number of modes, $P|_{max}$, we can compute off-line and store the total $N_{el}|_{max} - 1$ history matrices of size $(P|_{max}+1) \times (P|_{max}+1)$, which we can fetch later for any specific $N_{el} \leq N_{el}|_{max} $ and $P \leq P|_{max}$.

%******************************************************************************************
\begin{figure}[h]
	\center
	\includegraphics[width=0.5\linewidth]{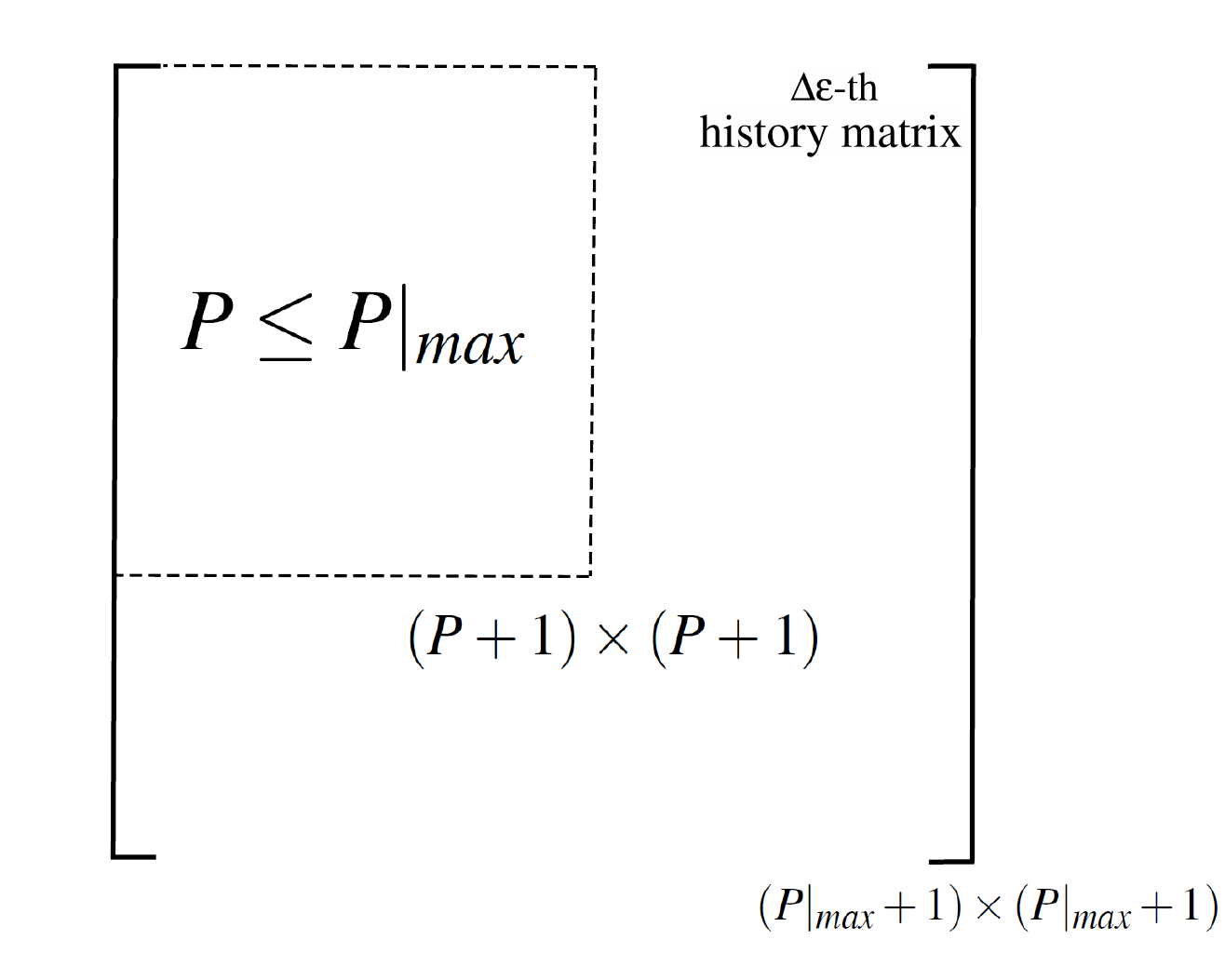}
	\caption{\label{Fig: History Sampling} History computation and retrieval.}
\end{figure}
% 

%\newpage

%%
%%%%%%%%%%%%%%%%%%%%%%%%%%%%%%%%%%%%%%%%%%
%\subsection{\textbf{H-Refinement for Stiff Problems}}
%\label{Sec: Stiff problems}
%%%%%%%%%%%%%%%%%%%%%%%%%%%%%%%%%%%%%%%%%%
%%
%There are different sources of singularity in the proposed problem. One is due to the force applied to the system, which can be captured by proper choice of basis and test functions. The second is due to the kernel of fractional derivative, which can cause singularities at the vicinity of boundaries. In this case, the singularity can be resolved by refining the grid within the boundary layer nearby the boundaries. In the next two subsections, we refine the grid using i) the kernel to obtain a non-uniform distribution of size of elements, and ii) a geometric series for the size of elements.
%
%

%
%%%%%%%%%%%%%%%%%%%%%%%%%%%%%%%%%%%%%%%%%
\subsection{\textbf{Non-Uniform Kernel-Based Grids}}
\label{Sec: kernel deriven grid}
\begin{figure}[h]
	\centering
	\begin{subfigure}{0.7\textwidth}
		\centering
		\includegraphics[width=1\linewidth]{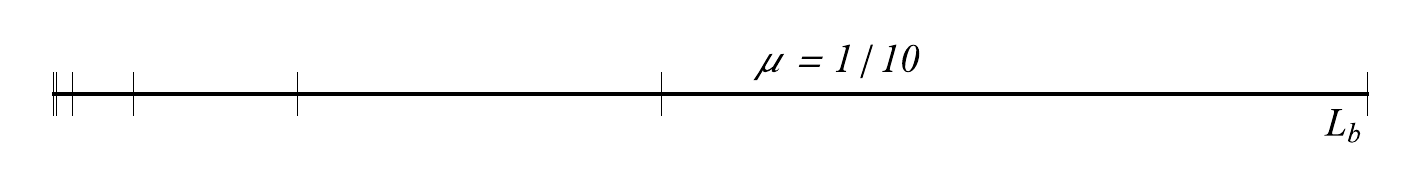}
		%		\caption{: Case I}
	\end{subfigure}
	\begin{subfigure}{0.7\textwidth}
		\centering
		\includegraphics[width=1\linewidth]{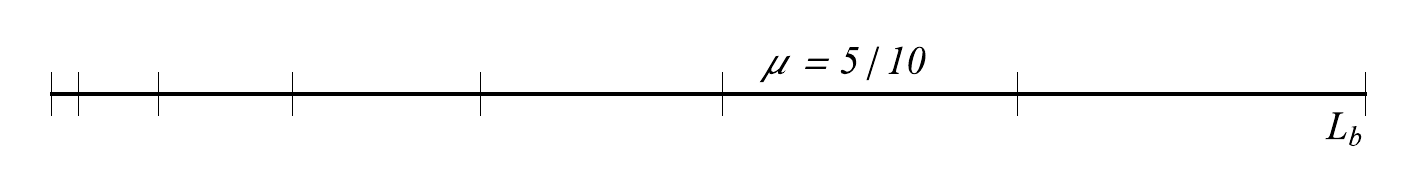}
		%		\caption{: Case III}
	\end{subfigure}
	\begin{subfigure}{0.7\textwidth}
		\centering
		\includegraphics[width=1\linewidth]{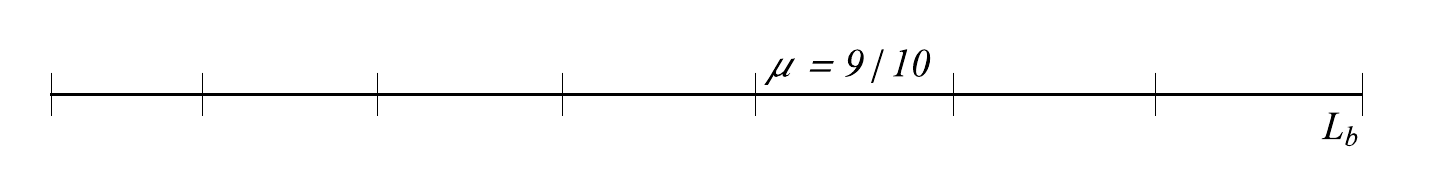}
		%		\caption{: Case III}
	\end{subfigure}
	\caption{\small{Non-uniform kernel-based grids for $N_b = 7$}.}
	\label{Fig: kernel driven grid}
\end{figure}

%%%%%%%%%%%%%%%%%%%%%%%%%%%%%%%%%%%%%%%%%
%
We present a non-uniform grid generation based on the power-law kernel in the definition of fractional derivative. There are different sources of singularity in the proposed problem that can be caused mainly due to the force function $f(x)$. However, even if the force term is smooth the underlying kernel of a fractional derivative leads to formation of singularities at the boundaries. Herein, we propose a new kernel-based grid generation method that considers a sufficiently small boundary layer at the vicinity of singular point and partitions that particular region non-uniformly. In this approach, we treat the kernel of the form $\frac{1}{x^{\sigma}}$ as a density function and then, we construct the grid such that the integral of kernel function over each element $\Omega_e \in [x_{e-1} , x_e]$ (in the boundary layer) is constant. 
%
%******************************************************************************************
\begin{figure}[h]
	\center
	\includegraphics[width=0.65\textwidth]{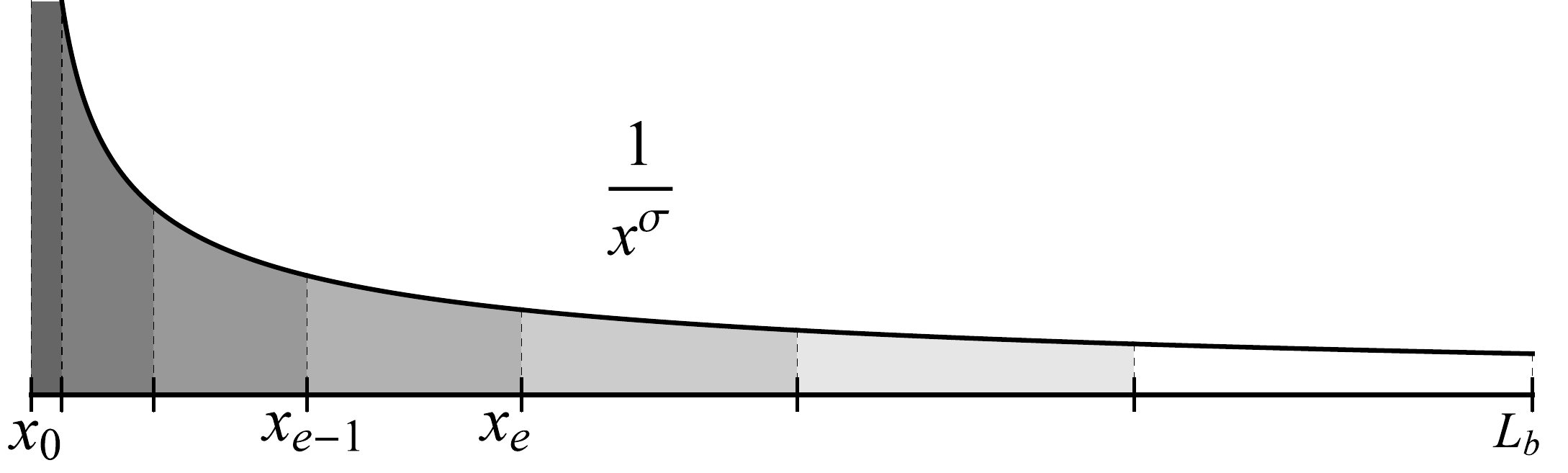}
	\caption{\label{Fig: nonuniform mesh} Kernel-based non-uniform grid in the boundary layer; $L_b$ and $N_b$ are the length of and the number of elements in the boundary layer, respectively.   }
\end{figure}
Since the operator is a left sided fractional derivative, we represent the non-uniform grid refinement at the left boundary. Let $L_b$ be the length of boundary layer and $\int_{0}^{L_b} \frac{1}{x^{\sigma}} \, dx = \frac{L_b^{1-\sigma}}{1-\sigma} = A$. Then, the integral over each element is
\begin{equation*}
	\frac{1}{A} \, \int_{x_{e-1}}^{x_e} \frac{1}{x^{\sigma}} \, dx = \frac{1}{L_b^{1-\sigma}} \, \left[ (x_{e-1} + \Delta x_e)^{1-\sigma} - x_{e-1}^{1-\sigma} \right] = \mathcal{C},
\end{equation*}
where $\Delta x_{e} = x_e - x_{e-1}$ and $\mathcal{C}$ is a constant. Thus,
\begin{equation*}
	\Delta x_{e} =  \left[ x_{e-1}^{1-\sigma} + \mathcal{C} \, L_b^{1-\sigma}  \right]^{\frac{1}{1-\sigma}}  - x_{e-1}.
\end{equation*}
Starting from $x_0 = 0$ and calculating the rest of grid locations successively, we obtain
\begin{align}
	\label{Eq: nonuniform mesh}
	x_e = \delta \,\, e^{\frac{1}{1-\sigma}} , \quad  e=1,2,\cdots,N_{b}, \quad \text{element numbers} ,
\end{align}
in which $\delta = L_b \,\, \mathcal{C}^{\frac{1}{1-\sigma}}$ and $N_b$ is the number of elements in the boundary layer. The constant $\mathcal{C}$ is obtained by the constraint $\sum_{e=1}^{N_{b}} \Delta x_e = L_b$ and hence,
\begin{align*}
	\mathcal{C} = \left(\sum_{e=1}^{N_{b}} \left[ e^{\frac{1}{1-\sigma}} - (e-1)^{\frac{1}{1-\sigma}} \right]  \right)^{\sigma - 1}.
\end{align*}
We consider $\sigma = 1 - \mu$ and thus when $\mu = 1$, we recover the uniform grid $x_e = \frac{L_b}{N_{b}} \,e$, where the kernel is $1$, $\mathcal{C} = \frac{1}{N_{b}}$, $\delta = \frac{L_b}{N_{b}}$. Fig. \ref{Fig: kernel driven grid} shows how the singularity in the kernel changes the non-uniformity in the grid.
%
%******************************************************************************************

\begin{figure}[h]
	\center
	\includegraphics[width=0.8\textwidth]{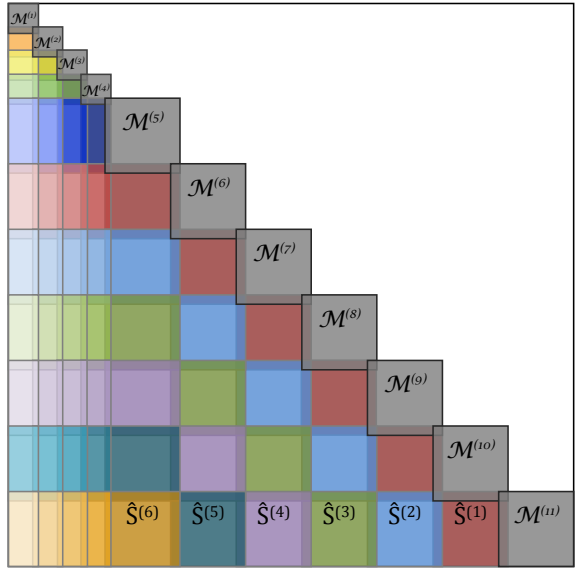}
	\caption{\label{Fig: Global Mat nonuniform}
		The assembled global matrix corresponding to $N_{el} = 11$ with $N_b = 4$ non-uniform boundary elements and $7$ uniform interior elements. In this global matrix, $\mathcal{M}^{(\varepsilon)} = \textbf{S}^{(\varepsilon)} - \lambda \textbf{M}^{(\varepsilon)}$, $\varepsilon =1,2,\cdots, N_{el}$, represents the local matrix, associated with the element $\Omega_{\varepsilon}$. The lower-triangle consists of three parts: 1) The small square $\frac{N_b(N_b-1)}{2}$ history matrices (interaction of boundary elements, $\varepsilon = 1,2,\cdots,N_b$). 2) The big square history matrices (interaction of interior elements, $\varepsilon = N_b+1,\cdots,N_{el}$). 3) The skinny rectangular $(N_{el}-N_b)N_b$ history matrices (interaction of boundary elements with interior elements).
	}
\end{figure}

We note that in the boundary layer, where the grid is non-uniform, equations \eqref{Eq: uniform mesh}-\eqref{Eq: h-I uniform} no longer hold. Thus, using \eqref{Eq: nonuniform mesh}, we obtain
%%
%\begin{align}
%%
%x_e + x_{e-1} & = \delta \,\, \left[ e^{\frac{1}{1-\mu}} + (e-1)^{\frac{1}{1-\mu}} \right] ,
%\\
%\nonumber
%x_e - x_{e-1} & = \delta \,\, \left[ e^{\frac{1}{1-\mu}} - (e-1)^{\frac{1}{1-\mu}} \right] , \quad
%%
%\end{align}
%
%
\begin{align}
	\label{Eq: History function general case nonuniform mesh}
	H_k^{(\varepsilon , e) }(\xi) 
	=
	\frac{-\mu}{\Gamma(1-\mu)} \frac{\delta}{2}   \left( \varepsilon^{\frac{1}{\mu}} - (\varepsilon-1)^{\frac{1}{\mu}} \right)
	\int_{-1}^{1}  \frac{\prescript{(2)}{}{\mathcal{P}}_{k+1}^{\mu}(\zeta) d\zeta}{\mathcal{Z}},
\end{align}
in which,
\begin{align*}
	\mathcal{Z}
	= \left( \frac{\delta}{2} \right)^{1+\mu}
	\Big[
	&
	\left( \varepsilon^{\frac{1}{\mu}} + (\varepsilon-1)^{\frac{1}{\mu}} \right) - \left( e^{\frac{1}{\mu}} + (e-1)^{\frac{1}{\mu}} \right)   
	\\
	&
	+ \left( \varepsilon^{\frac{1}{\mu}} - (\varepsilon-1)^{\frac{1}{\mu}} \right)  \zeta - \left( e^{\frac{1}{\mu}} - (e-1)^{\frac{1}{\mu}} \right)  \xi 
	\Big]^{1+\mu}.
\end{align*}
Therefore, by \eqref{Eq: global history mat} and \eqref{Eq: basis derivative}, the entries of the history matrix for the boundary layer elements, where $\varepsilon=1,2,\cdots,N_b$ and $e=1,2,\cdots,\varepsilon-1$, can be numerically obtained as 
\begin{equation}
	\label{Eq: global history mat-standard domain nonuniform}
	\boxed{
		\hat{\textbf{S}}^{(\varepsilon , e)}_{kp} =
		\int_{-1}^{1}
		\frac{d\psi_p }{d\xi} \, H_k^{(\varepsilon , e) }(\xi) \,\,   d\xi, \quad k,p =0,1,\cdots, P,
	}
\end{equation}
These matrices are the small squares in the upper left corner of Fig. \ref{Fig: Global Mat nonuniform} (interaction of boundary layer elements $e$ and $\epsilon$). For the interior elements, $\varepsilon=N_b+1, N_b+2, \cdots , N_{el}$, when $N_b+1 \leq e \leq \epsilon-1$, the grid is uniform and therefore, we use \eqref{Eq: global history mat-standard domain} to obtain the history matrices. These matrices are shown as the big squares in Fig. \ref{Fig: Global Mat nonuniform} (interaction of interior elements $e$ and $\epsilon$). However, when $1 \leq e \leq N_b$, the grid is non-uniform and we use \eqref{Eq: global history mat-standard domain nonuniform} to obtain the history matrices. These matrices are shown as skinny rectangles in Fig. \ref{Fig: Global Mat nonuniform} (interaction of interior elements with boundary layer elements).

%
%******************************************************************************************

In uniform grid generation, the history function \eqref{Eq: History function unifrom mesh} only depends on element difference $\Delta \epsilon$, which leads to a fast and efficient construction of history matrices (see Remark \ref{Rem: on the uniform mesh and history matrix}). However, in non-uniform kernel-based grid generation, this is not the case anymore and construction of history matrices is computationally expensive. Improving the history construction on non-uniform grids requires further investigations, to be done in our future works. 

%
%%%%%%%%%%%%%%%%%%%%%%%%%%%%%%%%%%%%%%%%%
\subsection{\textbf{Non-Uniform Geometrically Progressive Grids}}
\label{Sec: grid geometric series}
\begin{figure}[h]
	\center
	\includegraphics[width=1\textwidth]{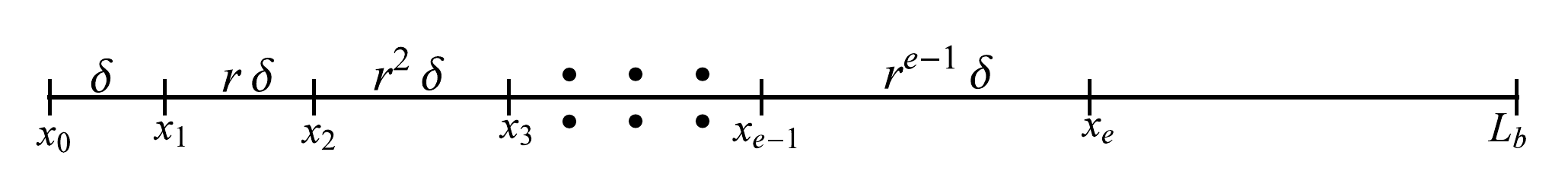}
	\caption{\label{Fig: nonuniform mesh geometric series} Non-uniform geometrically progressive grid.   }
\end{figure}

%%%%%%%%%%%%%%%%%%%%%%%%%%%%%%%%%%%%%%%%%
%
In addition to the non-uniform grid generation based on the kernel of fractional derivative, we consider a non-uniform grid using geometrically progressive series \cite{babuvska1994p,ainsworth2011posteriori}. In this case, the length of elements are increased by a constant factor $r$ (see Fig. \ref{Fig: nonuniform mesh geometric series}).
%
%******************************************************************************************
By considering the length of first element to be $\delta$, we construct the grid as $x_0 = 0 , \,\, x_1 = \delta , \,\, x_2 = \delta (1 + r) , \,\, x_3 = \delta (1 + r + r^2)$ and so on. Hence,
\begin{eqnarray}
	\label{Eq: nonuniform grid geometric series}
	x_e = \delta \sum_{i=0}^{e-1} r^i = \delta \frac{r^e - 1}{r-1}  , \quad e=1,2,\cdots,N_{b}.
\end{eqnarray}
Choosing $r$ and $N_{b}$, the constant $\delta$ is obtained by the constraint $x_{N_{b}} = L_b$, which gives $\delta = L_b \frac{r-1}{r^{N_{b}}-1}$. Since the grid is non-uniform, equations \eqref{Eq: uniform mesh}-\eqref{Eq: h-I uniform} do not hold anymore. Thus, using \eqref{Eq: nonuniform grid geometric series}, we obtain 
\begin{align}
	\label{Eq: History function general case geometric series nonuniform mesh}
	H_k^{(\varepsilon , e) }(\xi) 
	=
	\frac{-\mu}{\Gamma(1-\mu)} (\frac{2}{\delta})^{\mu}  r^{\Delta \varepsilon - \mu(e-1)}
	\int_{-1}^{1}  
	\frac{\prescript{(2)}{}{\mathcal{P}}_{k+1}^{\mu}(\zeta) d\zeta}{ \left[ \frac{r+1}{r-1} (r^{\Delta \varepsilon} - 1)  +  (\zeta - r^{\Delta \varepsilon} \xi) \right]^{1+\mu} },
\end{align}
where $\Delta \varepsilon = \varepsilon - e > 0$, denotes the \textit{element difference} between the current element $\varepsilon$ and the $e$-th element. Using the same expansion as in \eqref{Eq: poly-fractonomial expansion}, we can write \eqref{Eq: History function general case geometric series nonuniform mesh} as
\begin{eqnarray}
	\label{Eq: history-func H uniform-3}
	H^{(\varepsilon , e)}_k(\xi) = 
	\frac{-\mu}{\Gamma(1-\mu)} 
	(\frac{2}{\delta})^{\mu}  r^{\Delta \varepsilon - \mu(e-1)}
	\sum_{m=0}^{k} C_{km} 
	\tilde{h}^{(\varepsilon , e)}_m(\xi),
\end{eqnarray}
where the (modal) memory mode
\begin{align}
	\label{Eq: modal memory mode nonuniform geometric series}
	\tilde{h}^{(\varepsilon , e)}_m(\xi) = 
	\int_{-1}^{1} 
	\frac{(1-\zeta)^{\mu + m}  }{[ \frac{r+1}{r-1} (r^{\Delta \varepsilon} - 1)  +  (\zeta - r^{\Delta \varepsilon} \xi) ]^{1+\mu}} d\zeta, \quad m=0,1, \cdots, k,
\end{align}
can be obtained analytically using hypergeometric functions. Therefore, by \eqref{Eq: history-func H uniform-3} and \eqref{Eq: basis derivative}, the entries of the history matrix can be efficiently computed using the Gauss quadrature in \eqref{Eq: global history mat-standard domain}. The construction of the assembled global linear system is the same as kernel-based grid generation approach. We note that similar to uniform grid, in the non-uniform grid generation using the geometrical progression, the history functions depend on the element difference $\Delta \varepsilon = \varepsilon - e$, leading to a fast and efficient construction of history matrices.

%\newpage

%
%%%%%%%%%%%%%%%%%%%%%%%%%%%%%%%%%%%%%%%%%
\section{Petrov-Galerkin Method with Global Test Functions}
\label{Sec: PG Global Test}
%%%%%%%%%%%%%%%%%%%%%%%%%%%%%%%%%%%%%%%%%
%
In this section, similar to the case of local test functions, by substituting \eqref{Eq: u element-expansion} into \eqref{Eq: a def disc} and considering the global test function, given in \eqref{Eq: global test}, we obtain:
\begin{eqnarray}
	\nonumber
	&&
	\sum_{e=1}^{N_{el}}  
	\Big(  \sum_{p=0}^{P} \hat{u}^{(e)}_p \,\, \frac{d\psi_p(x) }{dx} \,,\, \prescript{RL}{x}{\mathcal{D}}_{L}^{\mu} \,v_k^{\varepsilon} (x) \Big)_{\Omega_e}  
	-
	\lambda 
	\sum_{e=1}^{N_{el}} 
	\Big( 
	\sum_{p=0}^{P} \hat{u}^{(e)}_p \,\, \psi_p (x) \,,\, \,v_k^{\varepsilon} (x)
	\Big)_{\Omega_e}  
	\\ \label{Eq: a def disc subs global-1}
	&=&  \sum_{e=1}^{N_{el}}  \Big( f \,,\,v_k^{\varepsilon} (x) \Big)_{\Omega_e}, \quad \varepsilon=1,2,\cdots,N_{el}, \quad k=0,1,\cdots,P .
\end{eqnarray}
Since the test function vanishes only $\forall x \in \Omega_e \neq \Omega_{\varepsilon}$ and $e>\varepsilon$, \eqref{Eq: a def disc subs global-1} reduces to
\begin{align}
	\nonumber
	&\sum_{e=1}^{\varepsilon} \sum_{p=0}^{P} \hat{u}^{(e)}_p 
	\Big(\frac{d\psi_p }{dx} \,,\, \prescript{RL}{x}{\mathcal{D}}_{x_{\varepsilon}}^{\mu} \,v_k^{\varepsilon} (x) \Big)_{\Omega_e}  -
	\lambda 
	\sum_{e=1}^{\varepsilon } \sum_{p=0}^{P} \hat{u}^{(e)}_p 
	\Big( \psi_p ( x  ) \,,\, \,
	v_k^{\varepsilon} (x) \Big)_{\Omega_e}  
	\\ \nonumber 
	&
	=
	\sum_{e=1}^{ \varepsilon } 
	\Big( f \,,\,v_k^{\varepsilon} (x) \Big)_{\Omega_e} .
\end{align}
By substituting \eqref{Eq: global test}, we obtain 
\begin{align}
	\nonumber
	&\sum_{e=1}^{\varepsilon} \sum_{p=0}^{P} \hat{u}^{(e)}_p 
	\Big(\frac{d\psi_p }{dx} \,,\, \prescript{RL}{x}{\mathcal{D}}_{x_{\varepsilon}}^{\mu} \,
	\prescript{(2)}{}{\mathcal{P}}_{k+1}^{\mu}( x^{1\sim \varepsilon}  )
	\Big)_{\Omega_e}  -
	\lambda 
	\sum_{e=1}^{\varepsilon } \sum_{p=0}^{P} \hat{u}^{(e)}_p 
	\Big( \psi_p ( x  ) \,,\, \,\prescript{(2)}{}{\mathcal{P}}_{k+1}^{\mu}( x^{1\sim \varepsilon}) \Big)_{\Omega_e}  
	\\ \nonumber 
	& =   
	\int_{0}^{x_{\varepsilon}}
	f(x)  \,\,\prescript{(2)}{}{\mathcal{P}}_{k+1}^{\mu}( x^{1\sim \varepsilon}  ) 
	dx, \quad \varepsilon=1,2,\cdots,N_{el}, \quad k=0,1,\cdots,P,
\end{align}
which can be written in the matrix form as
\begin{eqnarray}
	\label{Eq: global Matrix Form}
	\sum_{e=1}^{\varepsilon } \sum_{p=0}^{P} \hat{u}^{(e)}_p 
	\,\,
	\Big[
	\hat{\textbf{S}}^{(\varepsilon , e)}_{kp}
	- 
	\lambda 
	\,\,
	\hat{\textbf{M}}^{(\varepsilon , e)}_{kp}
	\Big] = \textbf{f}^{(\varepsilon)}_k, 
	\quad \varepsilon=1,2,\cdots,N_{el}, \quad k=0,1,\cdots,P,
\end{eqnarray}
where
\begin{eqnarray}
	\label{Eq: global Matrix def S}
	\hat{\textbf{S}}^{(\varepsilon , e)}_{kp} &=&
	\,\,
	\Big(\frac{d\psi^{\varepsilon}_p }{dx} \,,\,  \prescript{RL}{x}{\mathcal{D}}_{x_{\varepsilon}}^{\mu} \,
	\Big[ 
	\prescript{(2)}{}{\mathcal{P}}_{k+1}^{\mu}( x^{1\sim \varepsilon}  ) \, \Big]  \Big)_{\Omega_{\varepsilon}}, 
	\\ 
	\label{Eq: global Matrix def M}
	\hat{\textbf{M}}^{(\varepsilon , e)}_{kp} &=& \Big( \psi^{\varepsilon}_p ( x  ) \,,\, \,\prescript{(2)}{}{\mathcal{P}}_{k+1}^{\mu}( x^{1\sim \varepsilon}) \Big)_{\Omega_{\varepsilon}},  
	\\ 
	\label{Eq: global Matrix def f}
	\textbf{f}^{(\varepsilon)}_k &=& \int_{0}^{x_{\varepsilon}} f \,\,\prescript{(2)}{}{\mathcal{P}}_{k+1}^{\mu}( x^{1\sim \varepsilon}  ) dx.
\end{eqnarray}

\vspace{0.25cm}
\begin{rem}
	\label{Rem: on the global test function choice1}
	The benefit of choosing such global test functions is now clear since we can analytically evaluate $\prescript{RL}{x}{\mathcal{D}}_{x_{\varepsilon}}^{\mu} \,
	\prescript{(2)}{}{\mathcal{P}}_{k+1}^{\mu}( x^{1\sim \varepsilon}  )$. However, we note that this choice of test functions introduces ``extra'' work associated with the construction of the ``history mass matrix'' $\hat{\textbf{M}}^{(\varepsilon , e)}$, $\forall e=1,2,\cdots, \varepsilon-1$, when $\lambda \neq 0$. 
	
	%Another drawback of this choice of test functions would be the extra cost of carrying out the quadrature over the increasing-in-length domains of integration when calculating the load-vector $\textbf{f}^{(\varepsilon)}$. 
	%
\end{rem}
\vspace{0.25cm}
\begin{rem}
	\label{Rem: on the global test function choice2}
	The choice of global test functions leads to extra cost of quadrature carried out over the increasing-in-length domains of integration in \eqref{Eq: global Matrix def f}. Depending on the behaviour of the force-term $f(x)$, this approach might require adaptive/multi-element quadrature rules to obtain the corresponding entries of the desired precision. 
\end{rem}
%
%\vspace{0.25cm}
%

%
%%%%%%%%%%%%%%%%%%%%%%%%%%%%%%%%%%%%%%%%%
\subsection{\textbf{Elemental (Local) Operations}: the construction of $\textbf{f}^{(\varepsilon)}$}
\label{Sec: lcoal matrices for global test}
%%%%%%%%%%%%%%%%%%%%%%%%%%%%%%%%%%%%%%%%%
%
Here, the construction of the load-vector is the only operation that could be regarded as ``local operations''. Hence,  

\begin{eqnarray}
	\nonumber
	\textbf{f}^{(\varepsilon)}_k &=& 
	\int_{0}^{x_{\varepsilon }}
	f(x) \,\, \prescript{(2)}{}{\mathcal{P}}_{k+1}^{\mu}( x^{1\sim \varepsilon}  ) 
	dx
	\\ \nonumber
	&=& 
	(\frac{x_{\varepsilon} }{2})
	\int_{-1}^{1}
	(1-\zeta)^{\mu}
	f(\,  x^{1\sim \varepsilon}(\zeta)  \,)\, P^{\mu , -\mu}_{k}(\zeta) 
	%\frac{dx^{\varepsilon}}{d\zeta} 
	d\zeta,
\end{eqnarray}
and thus,
\begin{equation}
	\nonumber
	\boxed{
		\textbf{f}^{(\varepsilon)}_k  =
		(\frac{x_{\varepsilon}  }{2})
		\sum_{q=1}^{Q} w^{ \mu , 0}_q  \,\,f(\,  x^{1\sim \varepsilon}(\zeta_q)  \,)\, P^{\mu , -\mu}_{k} (\, z^{\mu , 0}_q \,) ,
	}
\end{equation}
where $\{ w^{\mu , 0}_q  \,,\, z^{\mu , 0}_q  \}_{q=1}^{Q}$ are the Gauss-Lobatto-Jacobi weights and points corresponding to the parameters $\alpha = \mu$ and $\beta = 0$.

%
%%%%%%%%%%%%%%%%%%%%%%%%%%%%%%%%%%%%%%%%%
\subsection{\textbf{Global Operations}: the construction of $\hat{\textbf{S}}^{(\varepsilon , e)}$ and $\hat{\textbf{M}}^{(\varepsilon , e)}$}
\label{Sec: global matrices for global test}
%%%%%%%%%%%%%%%%%%%%%%%%%%%%%%%%%%%%%%%%%
%
The corresponding stiffness and mass matrices are global in nature. We obtain their entries using proper Gauss quadrature rules.

%
%%%%%%%%%%%%%%%%%%%%%%%%%%%%%%%%%%%%%%%%%
\subsection{\textbf{Assembling the Global System with Global Test Functions}}
\label{Sec: Assembling Gloabl Test}
\begin{figure}[t]
	\center
	\includegraphics[width=0.55\textwidth]{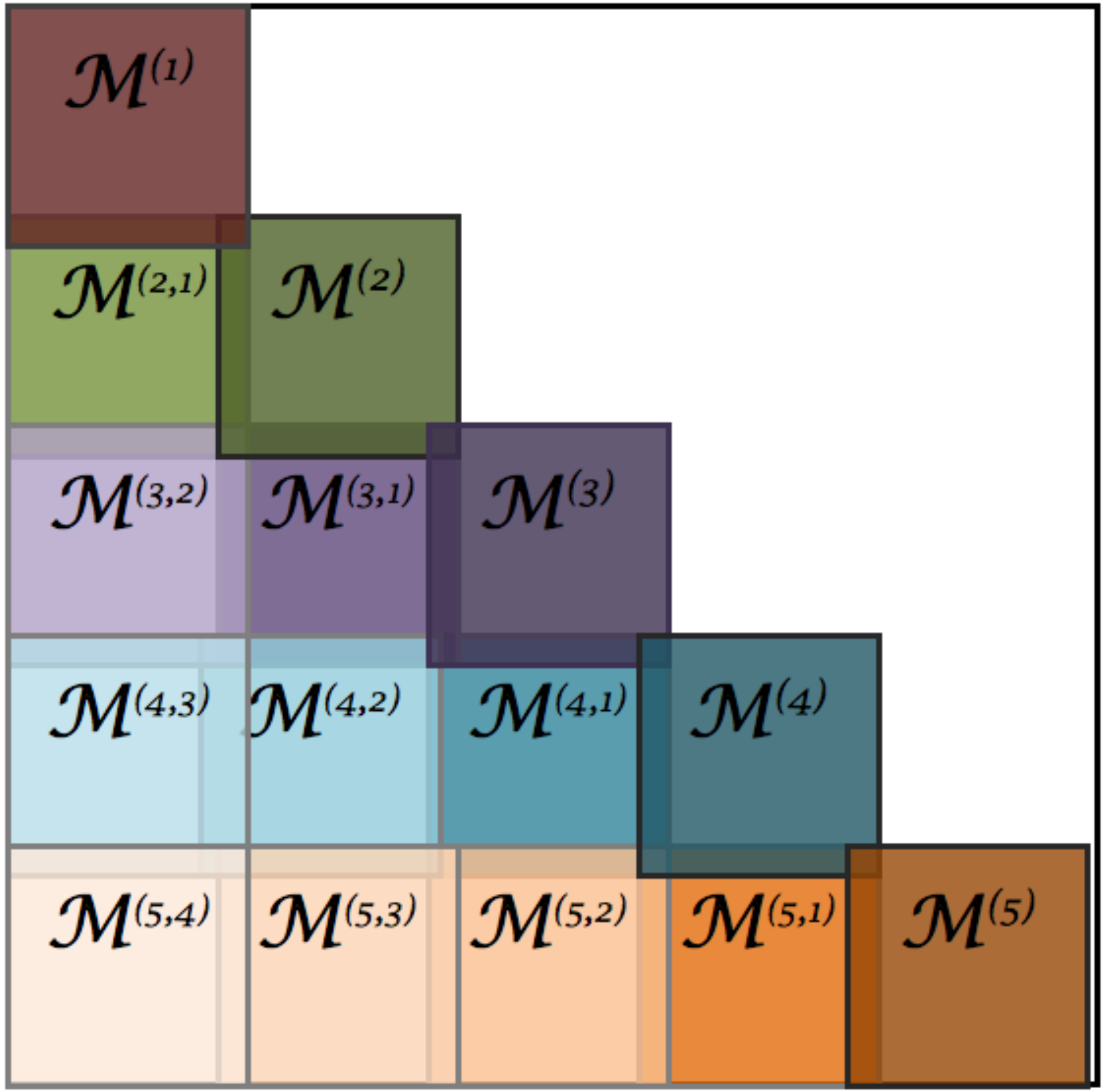}
	\caption{\label{Fig: Global Mat global test}
		The assembled global matrix corresponding to $N_{el} = 5$ elements when global test functions are employed. In this global matrix, $\mathcal{M}^{(\varepsilon)} = \hat{\textbf{S}}^{(\varepsilon)} - \lambda \hat{\textbf{M}}^{(\varepsilon)}$, $\varepsilon =1,2,\cdots, N_{el}$, represents the local matrix, associated with the element $\Omega_{\varepsilon}$. To fill the lower-triangular block matrices, we must construct $\frac{N_{el}(N_{el}-1)}{2}$ history matrices $\hat{ \textbf{S}   }^{( \varepsilon , e)}$.
	}
\end{figure}

%%%%%%%%%%%%%%%%%%%%%%%%%%%%%%%%%%%%%%%%%
%
We extend the notion of global linear system assembly by taking into account the presence of the history stiffness and mass matrices. We similarly impose the $C^0-continuity$ by employing the same ``mapping arrays'', $map[e][p]$, defined in \eqref{Eq: mapping array}. Let us define the $(P+1)\times(P+1)$ matrix
\begin{equation}
	\label{Eq: local linear system}
	\hat{\mathcal{M}}^{(\varepsilon , e)} = \hat{\textbf{S}}^{(\varepsilon , e)} - \lambda \hat{\textbf{M}}^{(\varepsilon , e)},
\end{equation}
$\forall \varepsilon, \,e$ fixed. Then, we assemble the corresponding global linear matrix $\mathbb{M}_G$ and the global load-vector $\mathbb{F}_G$ as follows:
\begin{mymathbox}
	\vspace{-0.25 in}
	\begin{eqnarray}
		\nonumber
		\textit{do}\,\,\varepsilon =1, N_{el} 
		\hspace{11cm} 
		\\ \nonumber
		\textit{do}\,\,k=1, P+1
		\hspace{10cm} 
		\\ \nonumber
		\mathbb{F}_G\Big[ \,\, map[\varepsilon][k]  \,\,   \Big] = \textbf{f}^{(\varepsilon)} [k]
		\hspace{8cm} 
		\\ \nonumber
		\textit{do}\,\,p=1, P+1
		\hspace{9.5cm} 
		\\ \nonumber
		%\mathbb{M}_G \Big[\,  map[\varepsilon][k]      \,\Big]\Big[\,   map[\varepsilon][p]  \,\Big]  = \mathbb{M}_G \Big[\,  map[\varepsilon][k]      \,\Big]\Big[\,   map[\varepsilon][p]  \,\Big] + \mathcal{M}^{(\varepsilon)}[k][p]
		%\hspace{-1cm} 
		%\\ \nonumber
		{\color{blue}
			\textit{do}\,\,e=1, {\color{red}\varepsilon}
			\hspace{9.50cm} 
		}
		\\ \nonumber
		{\color{blue}
			\mathbb{M}_G \Big[\,  map[  {\color{red}\varepsilon}  ][k]      \,\Big]\Big[\,   map[  {\color{red}e}  ][p]  \,\Big]  = \mathbb{M}_G \Big[\,  map[  {\color{red}\varepsilon}  ][k]      \,\Big]\Big[\,   map[  {\color{red}e}  ][p]  \,\Big]  + \hat{\mathcal{M}}^{(\varepsilon , e)}[k][p]
			\hspace{-1.5cm} 
		}
		%\\ \nonumber
		%{\color{blue}
		%End
		%}
		%\hspace{10.60 cm} 
		%\\ \nonumber
		%End
		%\hspace{11.2 cm} 
		%\\ \nonumber
		%End
		%\hspace{11.60 cm} 
		\\ \nonumber
		End
		\hspace{12.15 cm} 
	\end{eqnarray}
\end{mymathbox}
This leads to a linear system similar to that in \eqref{Eq: global linear system}, shown in Fig. \ref{Fig: Global Mat global test}, where the homogeneous Dirichlet boundary conditions are enforced in a similar fashion as before. We note that the scattering operation follows the same steps as explained in section \ref{Sec: Scatterign}.
%

%
%******************************************************************************************
% 

%%
%%%%%%%%%%%%%%%%%%%%%%%%%%%%%%%%%%%%%%%%%%
%\subsection{\textbf{Scattering Operation}}
%\label{Sec: Scatterign global test}
%%%%%%%%%%%%%%%%%%%%%%%%%%%%%%%%%%%%%%%%%%
%%
%It follows the same steps as explained in section \ref{Sec: Scatterign}.

%%%%%%%%%%%%%%%%%%%%%%%%%%%%%%%%%%%%%%%%%%%%%%%%%%%%%%%%%%%%%%%%%%%%%%%%%%%%%%%%%%%%%%%%%%%%%%%%%%%%%%%%%%%%%%%%%%%%%%%%%%%%%%%%%%%%%%%%%%%%%%%%%%%%%%%%%%%%%%%%%%%%%%%%%%%%%%%%%%%%%%%%%%%%%%%%%%%%%%%%%%%%%%%%%%%%%%%%%%%%%%%%%%%%%%%%%%%%%%%%%%%%%%%%%%%%%%%%%%%%%%%%%%%%%%%%%%%%%%%%%%%%%%%%%%%%%%%%%%%%%%%%%%%%%%%%

%\newpage
%
%%%%%%%%%%%%%%%%%%%%%%%%%%%%%%%%%%%%%%%%%
\section{\textbf{Numerical Examples}}
\label{Sec: NumExamp}
%%%%%%%%%%%%%%%%%%%%%%%%%%%%%%%%%%%%%%%%%
%
We consider numerical examples of the two PG schemes we have proposed. We provide examples of smooth and singular solutions with singularities at boundary points and in the interior domain, where we show the efficiency of developed schemes in capturing the singularities. We also perform the off-line computation of history matrices and show the improvement of computational cost. Moreover, we construct non-uniform kernel-based and geometrically progressive grids and present the success of the two approaches in capturing singular solutions. Furthermore, we investigate the non-local effects for different cases of history fading. In this section, we consider the computational domain $L=1$.

%
%%%%%%%%%%%%%%%%%%%%%%%%%%%%%%%%%%%%%%%%%
\subsection{\textbf{Smooth Problems}}
\label{Sec: smooth problem}
\begin{figure}[h]
	\centering
	\begin{subfigure}{0.41\textwidth}
		\centering
		\includegraphics[width=1\linewidth]{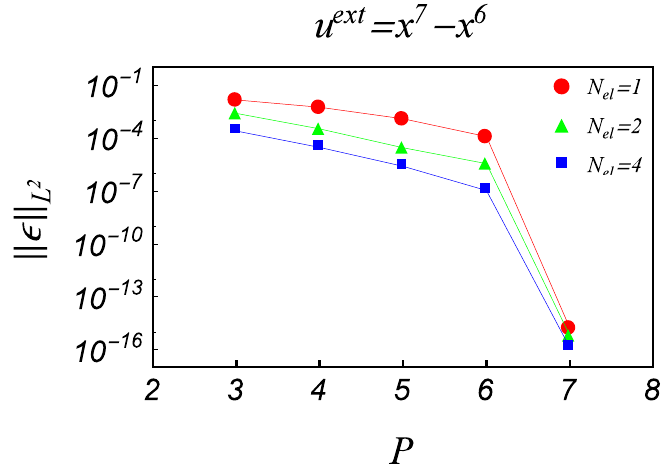}
		%		\caption{}
		%		\label{fig:smooth1}
	\end{subfigure}
	\begin{subfigure}{0.41\textwidth}
		\centering
		\includegraphics[width=1\linewidth]{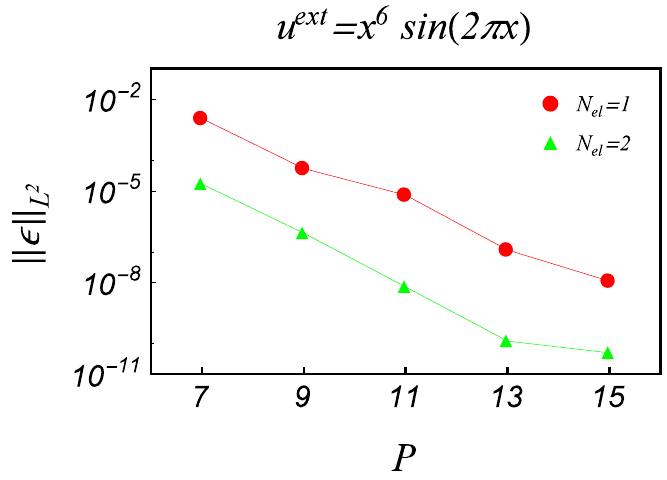}
		%		\caption{}
		%		\label{fig:smooth2}
	\end{subfigure}
	\caption{PG SEM with local basis/test functions. Plotted is the error with respect to the polynomial degree of each element (spectral order).}
	\label{Fig: smooth solutions}
\end{figure}

%%%%%%%%%%%%%%%%%%%%%%%%%%%%%%%%%%%%%%%%%
%
In the proposed schemes, the choice of bases functions are polynomials, enabling the scheme to accurately and efficiently approximate the smooth solutions over the whole domain. We consider two smooth solutions of the form $u^{ext} = x^7 - x^6$ and $u^{ext} = x^6 \sin(2\pi x)$. The corresponding force functions are obtained by substituting the exact solutions into \eqref{Eq: Helmholtz} (with $\lambda = 0$). By employing PG SEM, using local basis/test functions and local basis with global test functions (developed in Sec. \ref{Sec: PG Local Test} and Sec. \ref{Sec: PG Global Test}, respectively), we observe that the former leads to a better approximability and condition number. Fig. \ref{Fig: smooth solutions} presents the $L_2$-norm error of the PG SEM, employing local bases/test functions, where we show the exponential convergence of the scheme in approximating the two smooth solutions. The condition number of the resulting assembled global matrix, using the two developed schemes are also presented in Table \ref{Table: condition number LOCAL}. We show that the choice of local bases/test functions leads to a better conditioning for different number of elements and modes.

%
%******************************************************************************************
%
\begin{table}[h]
	\center
	\caption{\label{Table: condition number LOCAL} Condition number of the resulting assembled global matrix for the two choices of local bases/test functions (left) and local bases with global test functions (right) for different number of elements and modes.} 
	%
	%\vspace{0.1in}
	%
	\scalebox{0.7}{
		\begin{tabular}
			{l@{\hspace*{0.55 in}} c@{\hspace*{0.55 in}} c@{\hspace*{0in}}}
			\multicolumn{3}{c}{(Local Test Functions)}  \\[4pt]
			\hline
			\hline
			$P$ & $N_{el} = 2$ & $N_{el} = 10$   \\[4pt]
			\hline
			\\
			3 & 7.13   & 86.13
			\\
			\hline
			\\
			5 & 13.21  &  153.86
			\\
			\hline
			\\
			10 & 35.39    & 420.24
			\\
			\hline
			\hline
			%\vspace{0.2in}
			%
		\end{tabular}
	}
	%\vspace{0.2in}
	%
	\scalebox{0.7}{
		\begin{tabular}
			{l@{\hspace*{0.3 in}} c@{\hspace*{0.3 in}} c@{\hspace*{.0in}}}
			\multicolumn{3}{c}{(Global Test Function)} \\[4pt]
			\hline
			\hline
			$P$ & $N_{el} = 2$ & $N_{el} = 10$   \\[4pt]
			\hline
			\\
			3 &  3.46 $\times 10^{4}$    & 1.84 $\times 10^{16}$
			\\
			\hline
			\\
			5 &   4.3 $\times 10^{7}$    & 7.2 $\times 10^{16}$ 
			\\
			\hline
			\\
			10 &  2.73 $\times 10^{15}$  &  5.1 $\times 10^{17}$ 
			\\
			\hline
			\hline
			%\vspace{0.2in}
			%
		\end{tabular}
	}
\end{table}
\subsection{\textbf{History Retrieval}}
\label{Sec: History Retrieval}
%%%%%%%%%%%%%%%%%%%%%%%%%%%%%%%%%%%%%%%%%
%
As discussed in Sec. \ref{Sec: offline computation}, a large number of history matrices can be computed off-line, stored, and retrieved for later use. The retrieval process, compared to on-line construction of the history matrices, leads to higher computational efficiency. In this section, by considering $1000$ elements, we compute and store $999$ history matrices for different number of modes, $P=2,3,$ and $4$ (here $\mu = \frac{1}{2}$). Then, for different number of elements, we compute the CPU time required for constructing and solving the linear system, obtained by retrieving the stored history matrices from hard drive. We also compute the CPU time required for constructing and solving the linear system, obtained by on-line computation of the history matrices. Table \ref{Table: history retrieval} shows that in the case of $p=4$ and for $N_{el}=10$, $N_{el}=100$, $N_{el}=500$, and $N_{el}=1000$, the retrieval process is almost 4, 5, and 10 times faster, respectively. Thus, the higher $p$ is, the faster and more efficient the retrieval becomes.

\begin{table}[h]
	\center
	\caption{\label{Table: history retrieval} CPU time of constructing and solving the linear system based on off-line retrieval and on-line calculation of history matrices.} 
	%
	%	\vspace{0.2in}
	%
	\scalebox{0.6}{
		\begin{tabular}
			{ l c@{\hspace*{.30in}} c@{\hspace*{.30in}} c@{\hspace*{.30in}} c@{\hspace*{.30in}} c@{\hspace*{.30in}} c@{\hspace*{.30in}} c@{\hspace*{.30in}} c}
			\multicolumn{9}{c}{CPU Time} \\[4pt]
			\hline
			\hline
			& \multicolumn{2}{c}{$N_{el} = 10$}  &  \multicolumn{2}{c}{$N_{el} = 100$}  &  \multicolumn{2}{c}{$N_{el} = 500$} &  \multicolumn{2}{c}{$N_{el} = 1000$}\\[4pt]
			\hline
			\\
			$P$ & Off-line  & On-line   & Off-line  & On-line  & Off-line  & On-line  & Off-line  & On-line
			\\  
			& retrieval  & computation  & retrieval  & computation  & retrieval  & computation  & retrieval  & computation 
			\\
			\hline
			\\
			2 & 2.6520  & 7.2540  & 24.7105 & 83.5229  & 141.3525  & 429.3147  & 370.6895  & 790.3478
			\\
			\hline
			\\
			3 & 4.7580  & 18.9073  & 46.0826 & 161.8042  & 266.0441  & 1308.8327  & 746.4959  & 4423.7671
			\\
			\hline
			\\
			4 & 8.8140  & 32.2922  & 84.8645 & 499.9988  & 485.7715  & 5599.4062  & 1392.8705  & 14709.4902
			\\
			\hline
			\hline
			%			\vspace{0.2in}
			%
		\end{tabular}
	}
\end{table}
%
%%%%%%%%%%%%%%%%%%%%%%%%%%%%%%%%%%%%%%%%%%%%%%%%
%

%
%%%%%%%%%%%%%%%%%%%%%%%%%%%%%%%%%%%%%%%%%
\subsection{\textbf{Singular Problems}}
\label{Sec: singular problem}
%%%%%%%%%%%%%%%%%%%%%%%%%%%%%%%%%%%%%%%%%
%

The developed PG spectral element method, compared to single-domain spectral methods, further leads to accurate solutions even in the presence of singularities via \textit{hp}-refinements at the vicinity of singularities, while still employing smooth polynomial bases. The error in the boundary layer is controlled by considering sufficient number of modes in the boundary layer elements. The error in the interior domain is then improved by performing \textit{p}-refinement in those elements. In order to investigate the performance of the scheme in capturing a singularity, we consider three types of singularities, including: i) single-boundary singularity, ii) full-boundary singularity, and iii) interior singularity (when discontinuous force functions are applied). 

\begin{figure}[h]
	\center
	\includegraphics[width=0.66\textwidth]{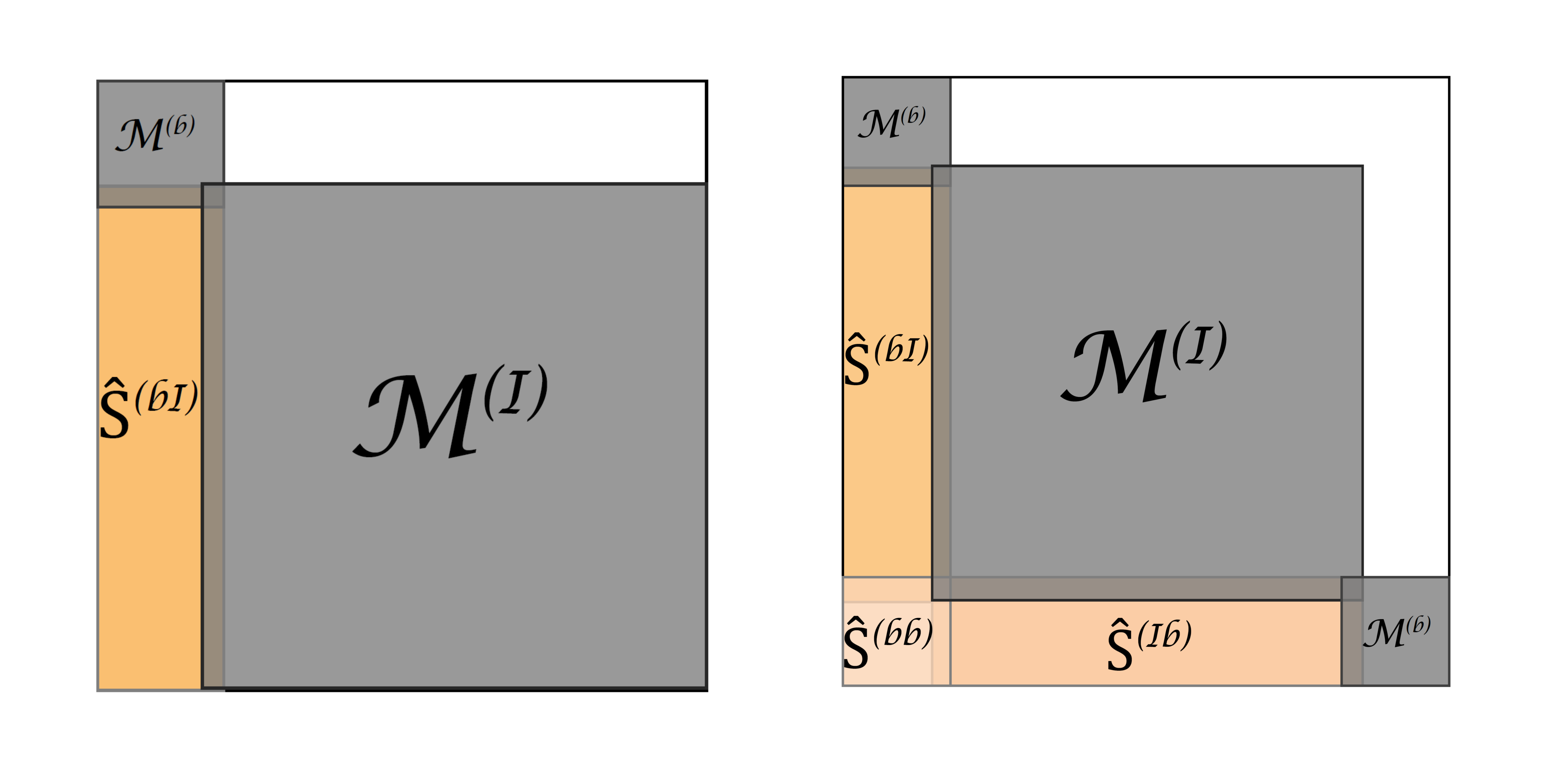}
	\caption{\label{Fig: Global Mat Boundary Layer - I} Schematic of global matrices corresponding to the case of singular solutions. (left): left boundary singularity, (right): left and right boundary singularities. $\hat{S}^{(bI)}$, $\hat{S}^{(Ib)}$, and $\hat{S}^{(bb)}$ denote the interaction of boundary/interior, interior/boundary and boundary/boundary elements, respectively.}
\end{figure}

\textbf{I) Single-Boundary Singularity:} we consider two singular solutions of the form $u^{ext} = (1-x) x^{2+\mu}$ and $u^{ext} = (1-x) x^{5+\mu}$ with left boundary singularity. We partition the domain into two non-overlapping elements, including one boundary element of length $L_b$ at the vicinity of singular point in addition to an interior element for the rest of computational domain. The schematic of corresponding global system is shown in Fig. \ref{Fig: Global Mat Boundary Layer - I} (left). Table \ref{Table: boundary layer - I} shows the exponential convergence of $L_2$-norm error in the interior domain. The error in the boundary layer element is then controlled by choosing sufficient number of modes in the boundary element. The results are obtained for the two cases of $L_b = 10^{-2} L$ and $L_b = 10^{-4} L$.  

\begin{table}[h]
	\center
	\caption{\label{Table: boundary layer - I} Single-Boundary Singularity: $L_2$-norm error in the boundary and interior elements using PG SEM with local basis/test functions. Here, $L_b$ represents the size of left boundary element, $P_b$ and $P_I$ denote the number of modes in the boundary and interior elements respectively. }
	\scalebox{0.7}{
		\begin{tabular}
			%
			%{r@{\hspace*{.30in}}c@{\hspace*{.30in}}c@{\hspace*{.30in}}c@{\hspace*{.30in}}c@{\hspace*{.30in}}c@{\hspace*{.30in}}c} 
			{l@{\hspace*{.30in}} c@{\hspace*{.30in}} c@{\hspace*{.30in}} c@{\hspace*{.30in}} }
			%			\hline
			%			\hline
			\\[-10 pt]
			\multicolumn{4}{c}{$u^{ext} = (1-x) x^{2+\mu}$, \,\, $\mu = 1/2$}  \\[4pt]
			\hline
			\hline
			\\[-10 pt]
			\multicolumn{4}{c}{\textbf{Boundary Element Error}}  \\[4pt]			
			\hline
			\\[-10 pt]
			$P_b$ & $L_b = 10^{-1}L$ & $L_b = 10^{-2}L$ & $L_b = 10^{-4}L$    \\[4pt]
			\hline
			\\[-10 pt]
			6  & $1.29387 \times 10^{-7}$  & $ 1.29634 \times 10^{-10}$  & $ 1.19525 \times 10^{-16}$  
			\\
			\hline
			\\[-10 pt]
			10 & $1.46601 \times 10^{-8}$ & $ 1.4193 \times 10^{-11}$  & $ 4.07955 \times 10^{-18}$ 
			\\
			\hline
			\hline
			\\[-10 pt]
			\multicolumn{4}{c}{\textbf{Interior Element Error, $P_b=10$}}  \\[4pt]
			\hline
			\\[-10 pt]
			$P_I$ & $L_b = 10^{-1}L$ & $L_b = 10^{-2}L$ & $L_b = 10^{-4}L$    \\[4pt]
			\hline
			\\[-10 pt]
			6 & $5.49133 \times 10^{-6}$ & $ 2.6893 \times 10^{-5}$  & $ 3.38957 \times 10^{-5}$  
			\\
			\hline
			\\[-10 pt]
			10 & $9.39045 \times 10^{-8}$ & $ 1.08594 \times 10^{-6}$  & $ 1.91087 \times 10^{-6}$  
			\\
			\hline
			\\[-10 pt]
			14 & $8.27224 \times 10^{-8}$ & $ 1.13249 \times 10^{-7}$  & $ 3.02065 \times 10^{-7}$ 
			\\
			\hline			  			  			  			  	 
			\hline
			\vspace{0.2in}
		\end{tabular}
	}
	\scalebox{0.7}{
		\begin{tabular}
			%
			%{r@{\hspace*{.30in}}c@{\hspace*{.30in}}c@{\hspace*{.30in}}c@{\hspace*{.30in}}c@{\hspace*{.30in}}c@{\hspace*{.30in}}c} 
			{l@{\hspace*{.30in}} c@{\hspace*{.30in}} c@{\hspace*{.30in}} c@{\hspace*{.30in}} }
			%			\hline
			%			\hline
			\\[-10 pt]
			\multicolumn{4}{c}{$u^{ext} = (1-x) x^{5+\mu}$, \,\, $\mu = 1/2$} \\[4pt]
			\hline
			\hline
			\\[-10 pt]
			\multicolumn{4}{c}{\textbf{Boundary Element Error}} \\[4pt]			
			\hline
			\\[-10 pt]
			$P_b$  & $L_b = 10^{-1}L$ & $L_b = 10^{-2}L$ &  $L_b = 10^{-4}L$   \\[4pt]
			\hline
			\\[-10 pt]
			6  & $3.94221 \times 10^{-11}$  & $ 2.96862 \times 10^{-17}$  & $ 4.8243 \times 10^{-29}$ 
			\\
			\hline
			\\[-10 pt]
			10 & $7.07024 \times 10^{-13}$  & $ 2.54089 \times 10^{-18}$  & $ 2.26939 \times 10^{-29}$ 
			\\
			\hline
			\hline
			\\[-10 pt]
			\multicolumn{4}{c}{\textbf{Interior Element Error, $P_b=10$}} \\[4pt]
			\hline
			\\[-10 pt]
			$P_I$ & $L_b = 10^{-1}L$  & $L_b = 10^{-2}L$ &  $L_b = 10^{-4}L$   \\[4pt]
			\hline
			\\[-10 pt]
			6  & $1.73622 \times 10^{-5}$  & $ 3.80264 \times 10^{-5}$  & $ 4.13249 \times 10^{-5}$ 
			\\
			\hline
			\\[-10 pt]
			10 & $1.3122 \times 10^{-9}$  & $ 8.76951 \times 10^{-9}$  & $ 1.10139 \times 10^{-8}$ 
			\\
			\hline
			\\[-10 pt]
			14 & $4.39611 \times 10^{-12}$  & $ 1.07775 \times 10^{-10}$  & $ 1.66044 \times 10^{-10}$ 
			\\
			\hline			  			  			  			  	 
			\hline
			%			\vspace{0.2in}
			%
		\end{tabular}
	}
\end{table}
%
%%%%%%%%%%%%%%%%%%%%%%%%%%%%%%%%%%%%%%%%%%%%%%%%%%
%

\textbf{II) Full-Boundary Singularity:} we consider the solution of the form $u^{ext} = (1-x)^{3+\mu_1} x^{3+\mu_2}$ with singular points at two ends, i.e. $x=0$ and $x=1$. Herein, we partition the domain into three non-overlapping elements including two boundary elements of length $L_b$ in the vicinity of singular points, and one interior element for the rest of domain. The schematic of corresponding global system is shown in Fig. \ref{Fig: Global Mat Boundary Layer - I} (right). Similar to previous example, the PG SEM can accurately capture the singularities at both ends, where increasing the number of modes in the interior element results in exponential convergence. Table \ref{Table: boundary layer - II} shows the $L_2$-norm error in the boundary layers and interior elements with two choices of $P_b=6 \,,\, 10$ and $L_b = 10^{-2}L \,,\, 10^{-4}L$.

\begin{table}[h]
	\center
	\caption{\label{Table: boundary layer - II} Full-Boundary Singularity: $L_2$-norm error in the boundary element (BE) and interior element (IE) by PG SEM with local basis/test functions. Here,  $u^{ext} = (1-x)^{3+\mu_1} x^{3+\mu_2}$ with $\mu_1=\frac{1}{4}, \,\, \mu_2=\frac{2}{3}$, $L_b$ represents the size of left and right boundary elements, $P_b$ and $P_I$ denote the number of modes in the boundary and interior elements respectively. }
	\scalebox{0.5}{
		\begin{tabular}
			%
			%{r@{\hspace*{.30in}}c@{\hspace*{.30in}}c@{\hspace*{.30in}}c@{\hspace*{.30in}}c@{\hspace*{.30in}}c@{\hspace*{.30in}}c} 
			{@{\hspace*{1.97in}}c@{\hspace*{1.97in}}}
			%			\hline
			%			\hline
			\\[-10 pt]
			$L_b = 10^{-2}L$ \qquad \qquad \quad \qquad \qquad \qquad \qquad \qquad \qquad \qquad \qquad \qquad \qquad $L_b = 10^{-4}L$   \\[4pt]
			%			 \multicolumn{7}{c}{$u^{ext} = (1-x)^{3+\mu_1} x^{3+\mu_2}$, $\quad$ $\mu_1=1/4 ,\,\, \mu_2=2/3$}   \\[4pt]
			\hline			  			  			  			  	 
			\hline
		\end{tabular}
	}
	\scalebox{0.5}{
		\begin{tabular}
			%
			%{r@{\hspace*{.30in}}c@{\hspace*{.30in}}c@{\hspace*{.30in}}c@{\hspace*{.30in}}c@{\hspace*{.30in}}c@{\hspace*{.30in}}c} 
			{l@{\hspace*{.30in}} c@{\hspace*{.30in}} c@{\hspace*{.30in}} c@{\hspace*{.30in}} }
			%			\hline
			%			\hline
			%			\\[-10 pt]
			%			\multicolumn{7}{c}{$u^{ext} = (1-x)^{3+\mu_1} x^{3+\mu_2}$, $\quad$ $\mu_1=1/4 ,\,\, \mu_2=2/3$}   \\[4pt]
			%			\hline
			%			\hline
			\\[-10 pt]
			\multicolumn{4}{c}{$P_b=6$}    \\[4pt]
			\hline
			\hline
			\\[-10 pt]
			$P_I$ & \textbf{Left BE Error} & \textbf{IE Error} & \textbf{Right BE Error}  \\[4pt]
			\hline
			\\[-10 pt]
			6 & $ 2.73893 \times 10^{-7}$ & $ 6.52605 \times 10^{-5}$ &  $ 3.51075 \times 10^{-6}$ 
			\\
			\hline
			\\[-10 pt]
			10 & $ 2.46964 \times 10^{-11}$  & $ 1.52215 \times 10^{-7}$ & $ 2.2902 \times 10^{-9}$
			\\
			\hline
			\\[-10 pt]
			14 & $ 3.08719 \times 10^{-12}$  & $ 9.30483  \times 10^{-9}$ & $ 2.69541 \times 10^{-10}$ 
			\\
			\hline
			\hline
			\\
			\multicolumn{4}{c}{$P_b=10$}   \\[4pt]
			\hline
			\hline
			\\[-10 pt]
			$P_I$ & \textbf{Left BE Error} & \textbf{IE Error} & \textbf{Right BE Error}  \\[4pt]
			\hline
			\\[-10 pt]
			6 & $ 2.73892 \times 10^{-7}$  & $ 6.52605 \times 10^{-5}$ &  $ 3.51075 \times 10^{-6}$ 
			\\
			\hline
			\\[-10 pt]
			10 & 2.48058 $ \times 10^{-11}$  & $ 1.52215 \times 10^{-7}$ & $ 2.29003 \times 10^{-9}$ 
			\\
			\hline
			\\[-10 pt]
			14 & $ 3.19684 \times 10^{-12}$  & $ 9.30511 \times 10^{-9}$ & $ 2.69506 \times 10^{-10}$ 
			\\
			\hline			  			  			  			  	 
			\hline
			%			\vspace{0.2in}
			%
		\end{tabular}
	}
	\scalebox{0.5}{
		\begin{tabular}
			%
			%{r@{\hspace*{.30in}}c@{\hspace*{.30in}}c@{\hspace*{.30in}}c@{\hspace*{.30in}}c@{\hspace*{.30in}}c@{\hspace*{.30in}}c} 
			{l@{\hspace*{.30in}} c@{\hspace*{.30in}} c@{\hspace*{.30in}} c@{\hspace*{.30in}} }
			%			\hline
			%			\hline
			%			\\[-10 pt]
			%			\multicolumn{7}{c}{$u^{ext} = (1-x)^{3+\mu_1} x^{3+\mu_2}$, $\quad$ $\mu_1=1/4 ,\,\, \mu_2=2/3$}   \\[4pt]
			%			\hline
			%			\hline
			\\[-10 pt]
			\multicolumn{4}{c}{$P_b=6$}   \\[4pt]
			\hline
			\hline
			\\[-10 pt]
			$P_I$ & \textbf{Left BE Error} & \textbf{IE Error} & \textbf{Right BE Error} \\[4pt]
			\hline
			\\[-10 pt]
			6 & $ 3.61679 \times 10^{-10}$  & $ 5.85397 \times 10^{-5}$ &  $ 4.60538 \times 10^{-8}$ 
			\\
			\hline
			\\[-10 pt]
			10 & $ 1.2676 \times 10^{-10}$  & $ 2.43295 \times 10^{-7}$ &  $ 1.62151 \times 10^{-10}$
			\\
			\hline
			\\[-10 pt]
			14 & $ 1.53993 \times 10^{-13}$  & $ 2.09933 \times 10^{-8}$ &  $1.97677 \times 10^{-11}$
			\\
			\hline
			\hline
			\\
			\multicolumn{4}{c}{$P_b=10$ }   \\[4pt]
			\hline
			\hline
			\\[-10 pt]
			$P_I$ & \textbf{Left BE Error} & \textbf{IE Error} & \textbf{Right BE Error} \\[4pt]
			\hline
			\\[-10 pt]
			6 & $3.61679 \times 10^{-10}$  & $ 5.85397 \times 10^{-5}$ &  $4.60538 \times 10^{-8}$
			\\
			\hline
			\\[-10 pt]
			10 & $ 1.2676 \times 10^{-12}$  & $ 2.43295 \times 10^{-7}$ &  $ 1.62151 \times 10^{-10}$
			\\
			\hline
			\\[-10 pt]
			14 & $ 1.53993 \times 10^{-13}$  & $ 2.09933 \times 10^{-8}$ &  $ 1.97677 \times 10^{-11}$
			\\
			\hline			  			  			  			  	 
			\hline
			%			\vspace{0.2in}
			%
		\end{tabular}
	}
\end{table}
%
%%%%%%%%%%%%%%%%%%%%%%%%%%%%%%%%%%%%%%%%%%%%%%%%%%
%

%
%******************************************************************************************
% 

%\newpage

\begin{figure}[t]
	\centering
	\begin{subfigure}{0.3\textwidth}
		\centering
		\includegraphics[width=1\linewidth]{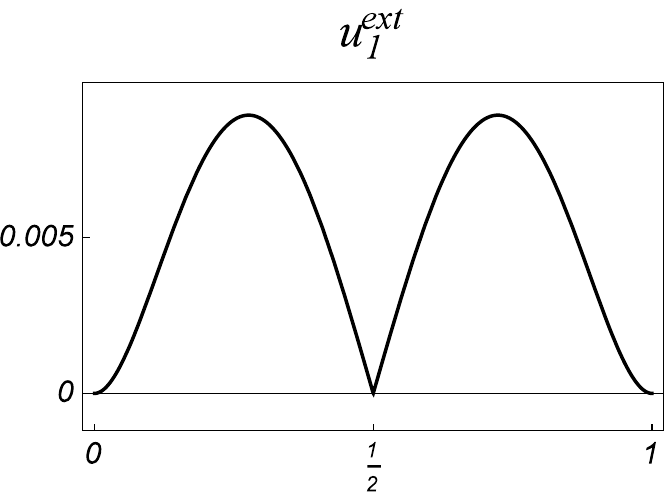}
		\caption{}
		\label{fig:sub1}
	\end{subfigure}
	\begin{subfigure}{0.28\textwidth}
		\centering
		\includegraphics[width=1\linewidth]{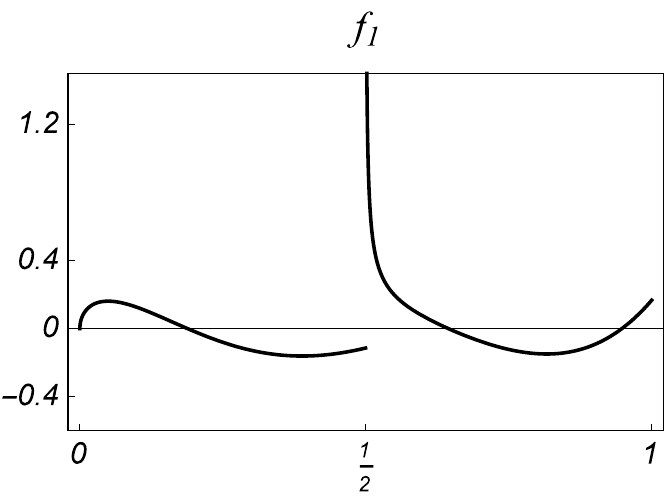}
		\caption{}
		\label{fig:sub2}
	\end{subfigure}
	\\
	\begin{subfigure}{0.3\textwidth}
		\centering
		\includegraphics[width=1\linewidth]{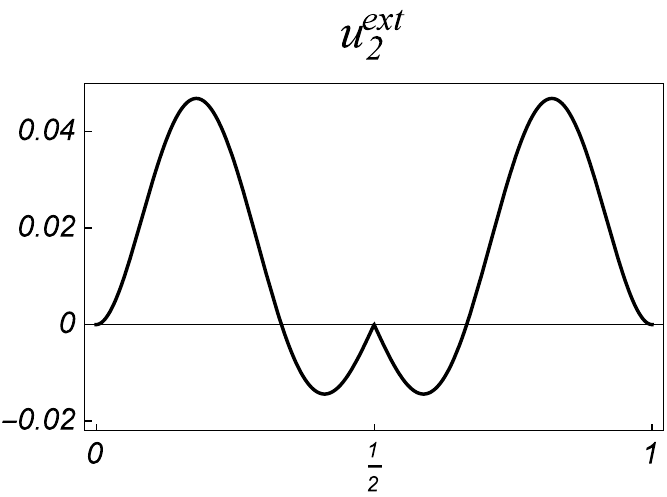}
		\caption{}
		\label{fig:sub3}
	\end{subfigure}
	\begin{subfigure}{0.28\textwidth}
		\centering
		\includegraphics[width=1\linewidth]{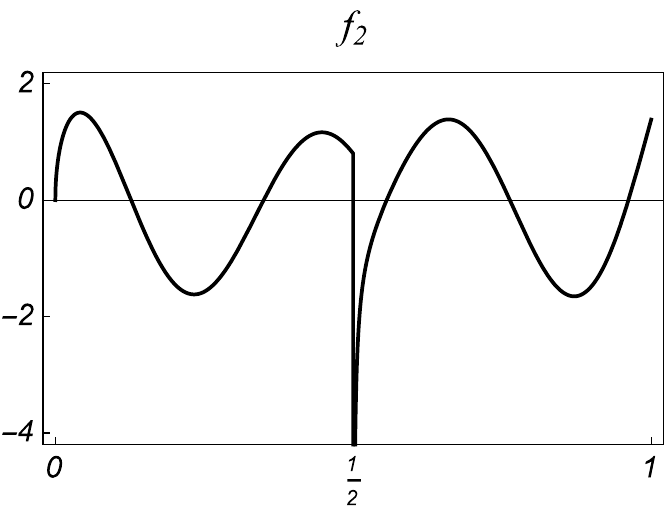}
		\caption{}
		\label{fig:sub4}
	\end{subfigure}
	\caption{Interior Singularity. (left): exact solutions, (right): the corresponding force functions. }
	\label{Fig: interior singularity}
\end{figure}

\begin{figure}[h]
	\centering
	\begin{subfigure}{0.35\textwidth}
		\centering
		\includegraphics[width=1\linewidth]{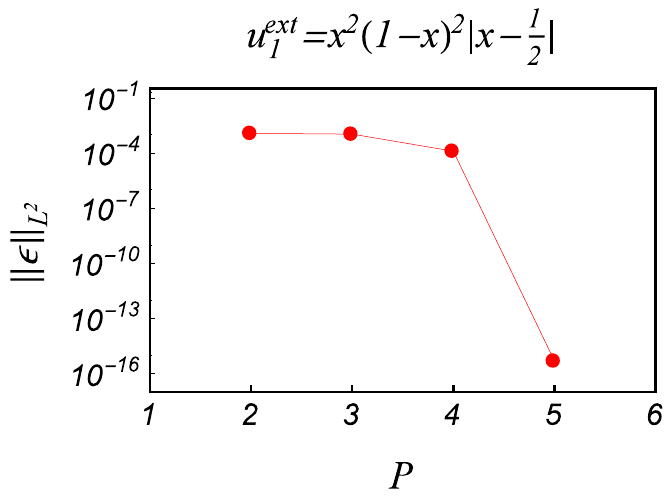}
		%		\caption{}
		%		\label{fig:midsingular1}
	\end{subfigure}
	\begin{subfigure}{0.35\textwidth}
		\centering
		\includegraphics[width=1\linewidth]{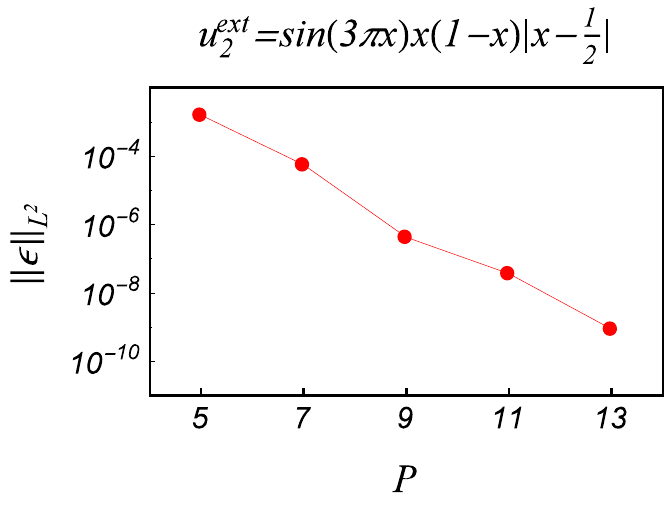}
		%		\caption{}
		%		\label{fig:midsingular2}
	\end{subfigure}	
	\caption{Interior Singularity: PG SEM with local basis/test functions. Plotted is the error with respect to spectral order in each element. }
	\label{Fig: interior singularity results}
\end{figure}

\textbf{III) Interior Singularity} (Discontinuous Force Function)\textbf{:} we consider the solutions with singularity in the middle of domain. The force function, obtained by substituting the solution into \eqref{Eq: Helmholtz}, is considered to be discontinuous at the point of singularity. Fig. \ref{Fig: interior singularity} shows the two exact solutions of the form $u_1^{ext} = x^2 \, (1-x)^2 \, |x-\frac{1}{2}|$ (top) and $u_2^{ext} = sin(3 \pi \, x) \, x \, (1-x) \, |x-\frac{1}{2}|$ (bottom) and their corresponding force functions. We partition the domain at the vicinity of singular point using two non-overlapping interior elements, in which the solution is smooth. The PG scheme with local basis/test functions is shown to be able to accurately capture the singularity in the middle of the domain. In the case of $u_1^{ext}$, we approximate the solution in the range of machine precision with $P=5$ within each element. We also show the exponential rate of convergence in the case of $u_2^{ext}$ by increasing the number of modes, $P$, in each element. The results are shown in Fig. \ref{Fig: interior singularity results}.

%
%******************************************************************************************
%

\newpage
%******************************************************************************************
%
%

%
%%
%\begin{table}[h]
%	\center
%	\caption{\label{Table: interior singularity} Interior Singularity: $L_2$-norm error by PG SEM with local basis/test functions.}
%	%
%	%	\vspace{0.2in}
%	%
%	\scalebox{0.6}{
%		%
%		\begin{tabular}
%			%
%			%{r@{\hspace*{.30in}}c@{\hspace*{.30in}}c@{\hspace*{.30in}}c@{\hspace*{.30in}}c@{\hspace*{.30in}}c@{\hspace*{.30in}}c} 
%			{l c}
%			\hline
%			\hline
%			\\[-10 pt]
%			\multicolumn{2}{c}{$u_1^{ext} = x^2 \, (1-x)^2 \, |x-\frac{1}{2}|$} \\[4pt]
%			\hline
%			\hline
%			\\[-10 pt]
%			$P$ & 5 \\[4pt]			
%			\hline
%			\\[-10 pt]
%			$L_2$-norm Error & $ 4.80753 \times 10^{-16}$ 
%			\\
%			\hline
%			\hline
%			\\
%			\\
%			%			\vspace{0.2in}
%			%
%		\end{tabular}
%		%
%	}
%	%
%	\scalebox{0.6}{
%		%
%		\begin{tabular}
%			%
%			%{r@{\hspace*{.30in}}c@{\hspace*{.30in}}c@{\hspace*{.30in}}c@{\hspace*{.30in}}c@{\hspace*{.30in}}c@{\hspace*{.30in}}c} 
%			{l@{\hspace*{.30in}} c@{\hspace*{.30in}} c@{\hspace*{.30in}} c@{\hspace*{.30in}} c@{\hspace*{.30in}} c}
%			\hline
%			\hline
%			\\[-10 pt]
%			\multicolumn{6}{c}{$u_2^{ext} = sin(3 \pi \, x) \, x \, (1-x) \, |x-\frac{1}{2}|$} \\[4pt]
%			\hline
%			\hline
%			\\[-10 pt]
%			$P$ & 5 & 7 & 9 & 11 & 13  \\[4pt]			
%			\hline
%			\\[-10 pt]
%			$L_2$-norm Error & $ 1.55714 \times 10^{-3}$ & $ 5.50791 \times 10^{-5}$  & $ 4.19702 \times 10^{-7}$ & $ 3.62245 \times 10^{-8}$ & $ 8.65968 \times 10^{-10}$
%			\\
%			\hline
%			\hline
%			%			\vspace{0.2in}
%			%
%		\end{tabular}
%		%
%	}
%	%
%\end{table}
%%
%%%%%%%%%%%%%%%%%%%%%%%%%%%%%%%%%%%%%%%%%%%%%%%%%%%
%%

%\newpage
%
%%%%%%%%%%%%%%%%%%%%%%%%%%%%%%%%%%%%%%%%%
\subsection{\textbf{Non-Uniform Grids}}
\label{Sec: stiff problem numerical example}
%%%%%%%%%%%%%%%%%%%%%%%%%%%%%%%%%%%%%%%%%
%

We consider a singular solution of the form $u^{ext} = (1-x) x^{1+\mu}$ (here $\mu = \frac{1}{10}$ and $\lambda = 0$) with singularity at the left boundary. In order to solve the problem, we consider three grid generation approaches with similar degrees of freedom, including one uniform and two non-uniform grids over the computational domain. The non-uniform grids are generated based on the power-law kernel in the definition of fractional derivative and the geometric progression series (discussed in Sec.\ref{Sec: kernel deriven grid} and Sec.\ref{Sec: grid geometric series}, respectively). Here, we choose $L_b = L$. Table \ref{Table: Nonuniform results} shows the $L_2$-norm error considering the uniform and non-uniform grids. We keep the total degrees of freedom fixed, but we increase the polynomial order $P$ in each simulation. The success of the non-uniform grid in providing more accurate results is observed, where fewer number of elements are used, while higher order polynomial are employed. We recall that the size of boundary layer has been set to its maximum possible length, i.e. $L_b = L$. Clearly, one can obtain even more accurate results when $L_b$ is set to much smaller length (e.g. $10^{-1}L$, $10^{-3}L$, \textit{etc.}). 

\begin{table}[h]
	\center
	\caption{\label{Table: Nonuniform results} $L_2$-norm error, using uniform and non-uniform grids. The exact singular solution is  $u^{ext} = (1-x) x^{1+\mu}$ with $\mu=1/10$.}
	%
	%	\vspace{0.2in}
	%
	\scalebox{0.6}{
		\begin{tabular}
			%
			%{r@{\hspace*{.30in}}c@{\hspace*{.30in}}c@{\hspace*{.30in}}c@{\hspace*{.30in}}c@{\hspace*{.30in}}c@{\hspace*{.30in}}c} 
			{l  c c c c }
			\hline
			\hline
			\\[-10 pt]
			&  Uniform Grid  &  Kernel-Based Non-Uniform Grid &  Geometrically Progressive Non-Uniform Grid\\[1pt]
			\hline
			\\[-10 pt]
			$N_{el} = 50$,  $P=2$ & $ 5.83943 \times 10^{-4}$ & $2.33461 \times 10^{-5}$ & $ 3.93956 \times 10^{-4}$ 
			\\[4pt]
			\hline
			\\[-10 pt]
			$N_{el} = 25$,  $P=4$ & $ 3.04739 \times 10^{-5}$ & $ 1.77458 \times 10^{-7}$ & $ 1.38755 \times 10^{-6}$ 
			\\[4pt]
			\hline
			\\[-10pt]
			$N_{el} = 10$, $P=10$ & $ 1.39586 \times 10^{-5}$ & $ 2.10813 \times 10^{-9}$ & $ 1.45695 \times 10^{-9}$ 
			\\
			\hline
			\hline
			%			\vspace{0.2in}
			%
		\end{tabular}
	}
\end{table}
%
%%%%%%%%%%%%%%%%%%%%%%%%%%%%%%%%%%%%%%%%%%%%%%%%%%
%

%%
%\begin{table}[h]
%	\center
%	\caption{\label{Table: Nonuniform results} $L_2$-norm error, using uniform and non-uniform grids. The exact solution is  $u^{ext} = (1-x) x^{1+\mu}$ with $\mu=1/10$.}
%	%
%	%	\vspace{0.2in}
%	%
%	\scalebox{0.6}{
%		%
%		\begin{tabular}
%			%
%			%{r@{\hspace*{.30in}}c@{\hspace*{.30in}}c@{\hspace*{.30in}}c@{\hspace*{.30in}}c@{\hspace*{.30in}}c@{\hspace*{.30in}}c} 
%			{l@{\hspace*{0.3in}} c c c c c }
%			\hline
%			\hline
%			\\[-10 pt]
%			 &  & Uniform Grid  &  Kernel Driven Non-Uniform Grid &  Geometrically Progressive Non-Uniform Grid\\[1pt]
%			\hline
%			\\[-10 pt]
%			\multirow{2}{1in}{$N_{el} = 25$ \\ $P=4$} & Projection error & $ 3.04739 \times 10^{-5}$ & $ 1.77458 \times 10^{-7}$ & $ 1.38755 \times 10^{-6}$ 
%			\\[4pt]
%			 & Solution error & $ 2.20784 \times 10^{-5}$ & $  \times 10^{-}$ & $ 1.04387 \times 10^{-3}$ 
%			\\[4pt]
%			\hline
%			\\[-10pt]
%			\multirow{2}{1in}{$N_{el} = 10$\\ $P=10$} & Projection error & $ 1.39586 \times 10^{-5}$ & $ 2.10813 \times 10^{-9}$ & $ 1.45695 \times 10^{-9}$ 
%			\\[4pt]
%			& Solution error & $1.56978  \times 10^{-5}$ & $  \times 10^{-}$ & $ 8.98276  \times 10^{-4}$ 
%			\\
%			\hline
%			\hline
%			%			\vspace{0.2in}
%			%
%		\end{tabular}
%		%
%	}
%	%
%\end{table}
%%
%%%%%%%%%%%%%%%%%%%%%%%%%%%%%%%%%%%%%%%%%%%%%%%%%%%
%%

%\newpage

%
%%%%%%%%%%%%%%%%%%%%%%%%%%%%%%%%%%%%%%%%%
\subsection{\textbf{A Systematic Memory Fading Analysis}}
\label{Sec: memory fading analysis}
%%%%%%%%%%%%%%%%%%%%%%%%%%%%%%%%%%%%%%%%%
%

In order to investigate the effect of truncating the history matrices, we perform a systematic memory fading analysis. 

In \textit{full} memory fading, we fade the memory by truncating the history matrices, i.e., we consider the full history matrices up to some specific number and then truncate the rest of history. For instance, we consider up to the first $4$ history matrices for each element and thus compute $\hat{\textbf{S}}^{1}$, $\hat{\textbf{S}}^{2}$, $\hat{\textbf{S}}^{3}$ and $\hat{\textbf{S}}^{4}$, and truncate the rest $N_{el}-1-4$ matrices; see Fig. \ref{Fig: Global Mat} for better visualization.

In \textit{partial} memory fading, we fade the memory by partially computing the history matrices. Similar to the \textit{full} memory fading, we consider the full history matrices up to some specific number, however, for the rest of history matrices we partially compute the entries of matrices. In partial memory fading, we consider three different cases as follows.

\begin{figure}[t]
	\centering
	\begin{subfigure}{0.15\textwidth}
		\centering
		\includegraphics[width=1\linewidth]{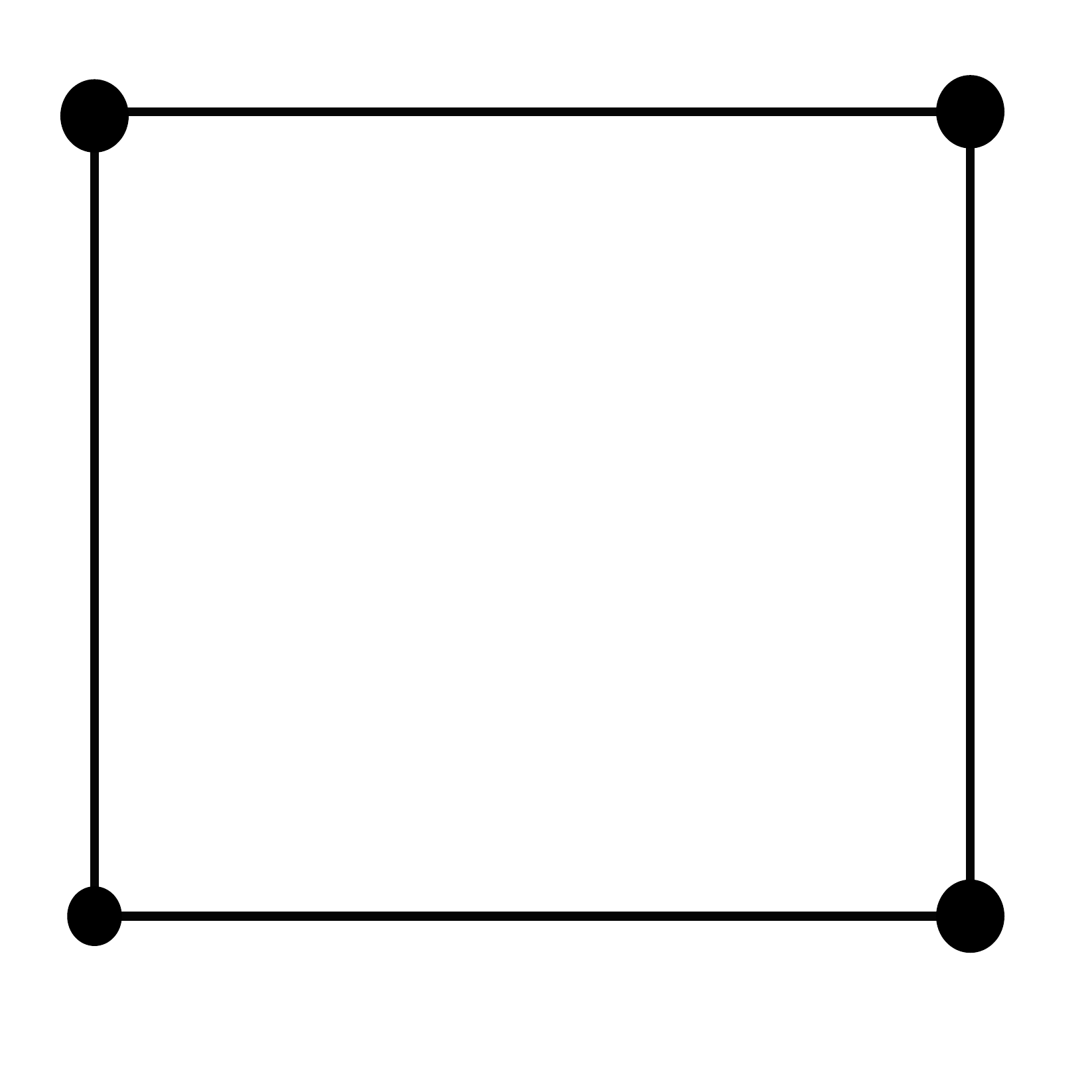}
		\caption{: Case I}
		\label{fig: Memory fading case I}
	\end{subfigure}
	\begin{subfigure}{0.15\textwidth}
		\centering
		\includegraphics[width=1\linewidth]{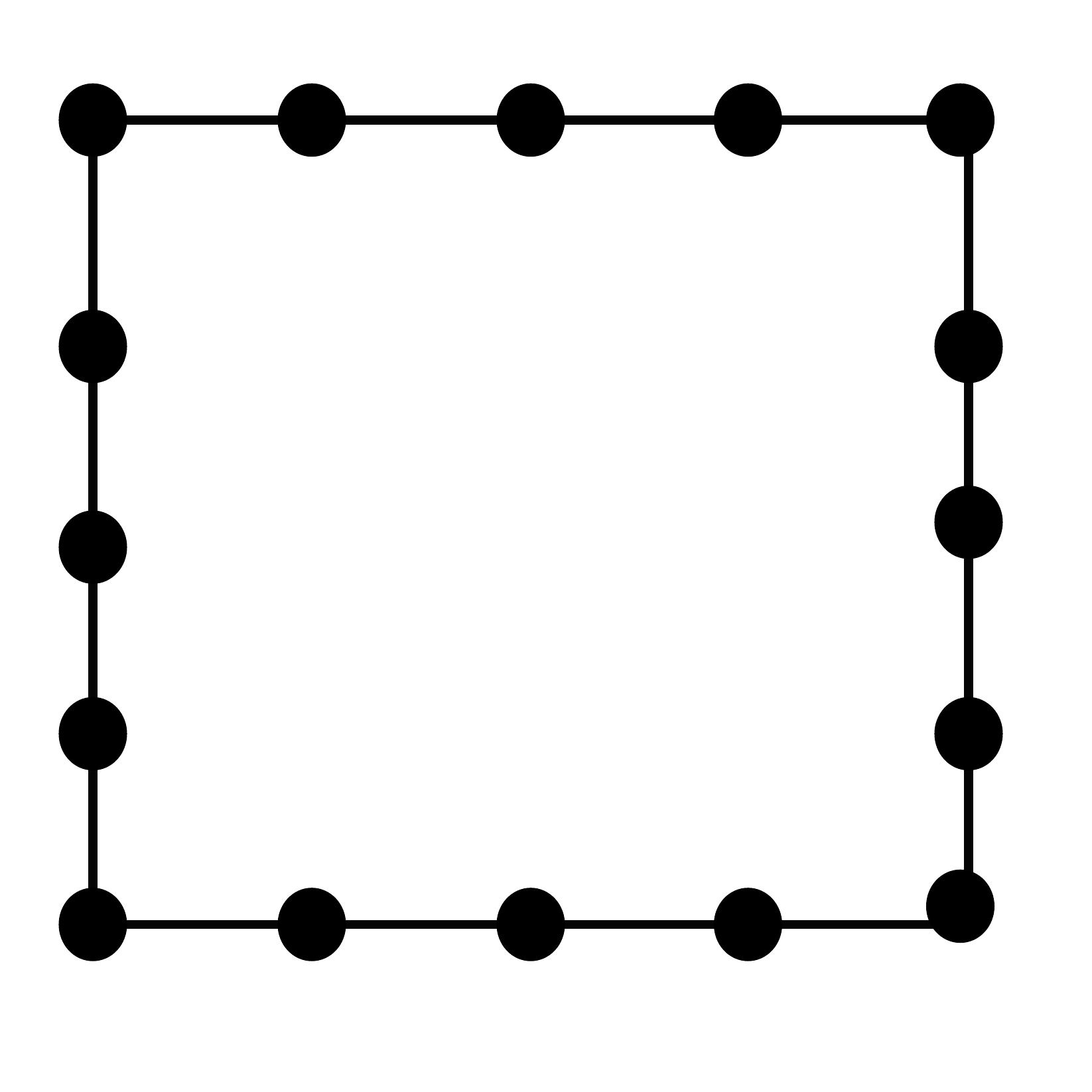}
		\caption{: Case II}
		\label{fig: Memory fading case II}
	\end{subfigure}
	\begin{subfigure}{0.15\textwidth}
		\centering
		\includegraphics[width=1\linewidth]{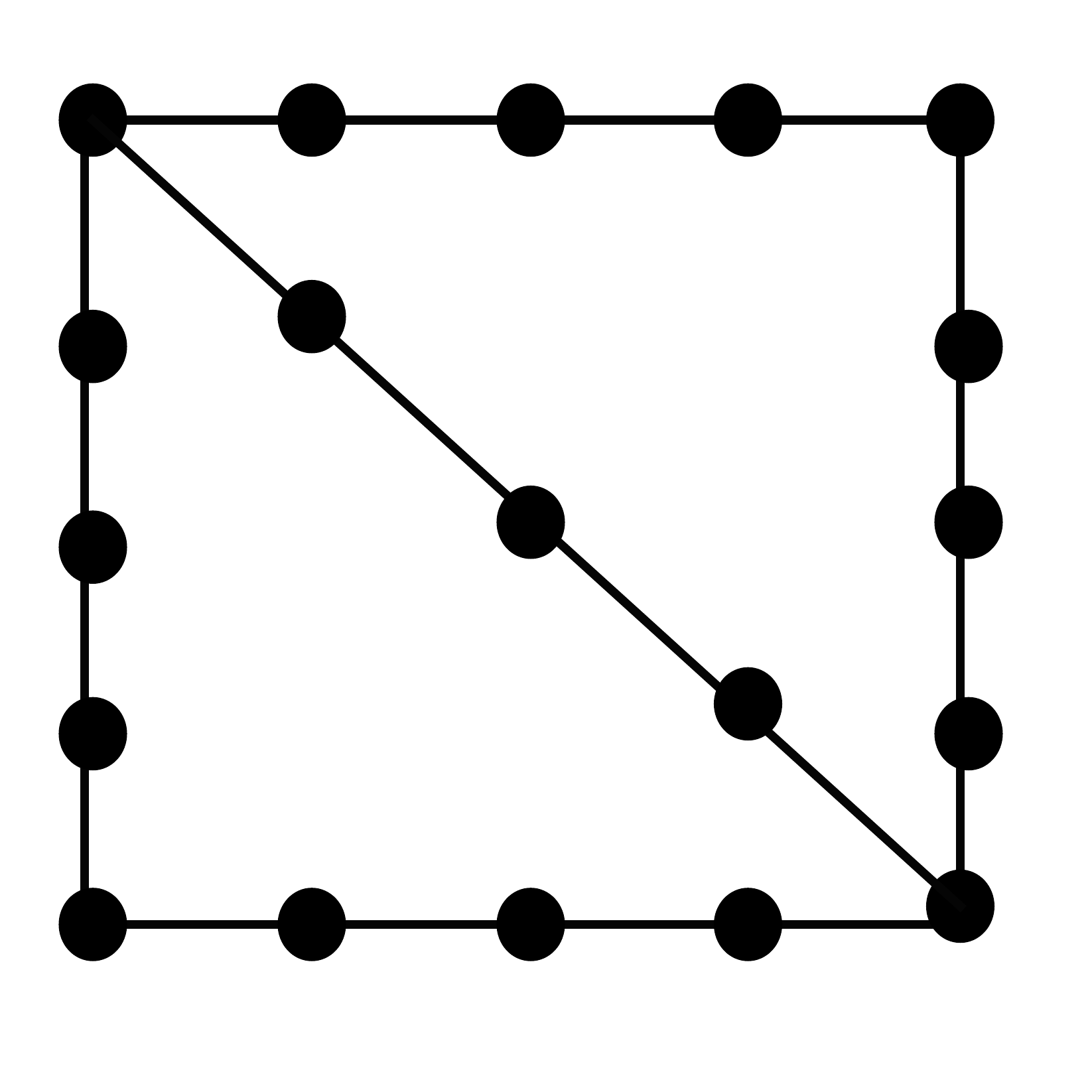}
		\caption{: Case III}
		\label{fig: Memory fading case III}
	\end{subfigure}
	\caption{\small{Memory fading: (a) B-B interaction, the corner entries (b) B-B and B-I interaction, the boundary entries (c) B-B, B-I and S-I interaction, boundary and diagonal entries}}
	\label{Fig:  Memory fading}
\end{figure}

\begin{itemize}
	\item \textbf{Case I:} Boundary-Boundary (B-B) interaction. In this case, we only consider the interactions of boundary mode and boundary test functions, i.e., $p = 0 , P $ and $k = 0 , P$, and thus, only compute the corner entries (See Fig. \ref{fig: Memory fading case I}).
	\item \textbf{Case II:} Boundary-Boundary (B-B) and Boundary-Interior (B-I) interaction. In addition to the corner entries, here we also consider the interaction of boundary mode/test functions with the interior test/mode functions, i.e., 
	\begin{align*}
		\begin{cases*}
			k=0, \quad p=0,1,\cdots,P, \quad \text{and} \quad  k=P, \quad p=0,1,\cdots,P
			\\
			p=0, \quad k=0,1,\cdots,P, \quad \text{and} \quad p=P, \quad k=0,1,\cdots,P,
		\end{cases*}
	\end{align*}
	and thus, we compute the boundary entries (See Fig. \ref{fig: Memory fading case II}).
	\item \textbf{Case III:} Boundary-Boundary (B-B), Boundary-Interior (B-I), Self-Interior (S-I) interaction. In addition to the last two cases, we consider the interaction of each mode with its corresponding test function and thus, we compute the boundaries as well as the diagonal entries (See Fig. \ref{fig: Memory fading case III}). 
\end{itemize}

%******************************************************************************************
%

%\newpage

Tables \ref{Table: full memory fading local} and \ref{Table: partial memory fading local} show the $L_2$-norm error for cases of full and partial memory fading. It is clear from the computed norms that even in the case of fading memory, we can still accurately obtain the approximation solution, however with a proportional loss of accuracy depending on the lack of modal interaction.

%
%%%%%%%%%%%%%%%%%%%%%%%%%%%%%%%%%%%%%%%%%%%%%%%%%%
%
\begin{table}[h]
	\center
	\caption{\label{Table: full memory fading local} Full history fading: $L_2$-norm error using PG SEM with local basis/test functions, where $u^{ext} = x^7 -x^6$, $N_{el} = 19$, $P=6$. The first column in the table shows the number of fully faded history matrices.} 
	%
	%	\vspace{0.2in}
	\scalebox{0.6}{
		\begin{tabular}
			%
			%{r@{\hspace*{.30in}}c@{\hspace*{.30in}}c@{\hspace*{.30in}}c@{\hspace*{.30in}}c@{\hspace*{.30in}}c@{\hspace*{.30in}}c} 
			{c c@{\hspace*{.30in}} c@{\hspace*{.30in}} c}
			\multicolumn{4}{c}{Full fading}\\
			\hline
			\hline
			\# faded history matrices & $\mu =1/10$ & $\mu =1/2$&  $\mu = 9/10$  \\[4pt]
			\hline
			0  & $ 9.26034 \times 10^{-12}$  & $ 2.31391 \times 10^{-11}$  & $ 4.24903 \times 10^{-9}$
			\\
			\hline
			2  & $ 7.8905 \times 10^{-11}$  & $ 1.26365 \times 10^{-10}$  & $ 4.25456 \times 10^{-9}$
			\\
			\hline
			5  & $ 1.42423 \times 10^{-8}$  & $ 6.39474 \times 10^{-8}$  & $ 1.95976 \times 10^{-8}$
			\\
			\hline
			8  & $ 2.69431 \times 10^{-7}$  & $ 2.47423 \times 10^{-6}$  & $ 8.45001 \times 10^{-7}$
			\\
			\hline
			11   & $ 2.09737 \times 10^{-6}$  & $ 3.19995 \times 10^{-5}$  & $ 1.37959 \times 10^{-5}$
			\\
			\hline
			14   & $ 9.07427 \times 10^{-6}$  & $ 2.44911 \times 10^{-4}$  & $ 1.40684 \times 10^{-4}$
			\\
			\hline			
			17   & $ 2.94001 \times 10^{-5}$  & $ 1.39043 \times 10^{-3}$  & $ 1.6001 \times 10^{-3}$
			\\
			\hline
			\hline
			\vspace{0.2in}
		\end{tabular}
	}
\end{table}
%
%%%%%%%%%%%%%%%%%%%%%%%%%%%%%%%%%%%%%%%%%%%%%%%%%%
%

%
%%%%%%%%%%%%%%%%%%%%%%%%%%%%%%%%%%%%%%%%%%%%%%%%%%
%
\begin{table}[h]
	\center
	\caption{\label{Table: partial memory fading local} Partial history fading: $L_2$-norm error using PG SEM with local basis/test functions, where $u^{ext} = x^7 -x^6$, $N_{el} = 19$, $P=6$. The first column in the tables shows number of partially faded history matrices.} 
	%
	%	\vspace{0.2in}
	%
	\scalebox{0.6}{
		\begin{tabular}
			%
			%{r@{\hspace*{.30in}}c@{\hspace*{.30in}}c@{\hspace*{.30in}}c@{\hspace*{.30in}}c@{\hspace*{.30in}}c@{\hspace*{.30in}}c} 
			{c c@{\hspace*{.30in}} c@{\hspace*{.30in}} c}
			\multicolumn{4}{c}{Partial fading case I}\\
			\hline
			\hline
			\# faded history matrices & $\mu =1/10$ & $\mu =1/2$&  $\mu = 9/10$  \\[4pt]
			\hline
			0  & $ 9.26034 \times 10^{-12}$  & $ 2.31391 \times 10^{-11}$  & $ 4.24903 \times 10^{-9}$
			\\
			\hline
			2  & $ 7.8905 \times 10^{-11}$  & $ 1.26365 \times 10^{-10}$  & $ 4.25456 \times 10^{-9}$
			\\
			\hline
			5  & $ 1.42423 \times 10^{-8}$  & $ 6.39474 \times 10^{-8}$  & $ 1.95976 \times 10^{-8}$
			\\
			\hline
			8  & $ 2.69431 \times 10^{-7}$  & $ 2.47423 \times 10^{-6}$  & $ 8.45001 \times 10^{-7}$
			\\
			\hline
			11   & $ 2.09737 \times 10^{-6}$  & $ 3.19995 \times 10^{-5}$  & $ 1.37959 \times 10^{-5}$
			\\
			\hline			
			14   & $ 9.07427 \times 10^{-6}$  & $ 2.44911 \times 10^{-4}$  & $ 1.40684 \times 10^{-4}$
			\\
			\hline
			17   & $ 2.94001 \times 10^{-5}$  & $ 1.39043 \times 10^{-3}$  & $ 1.6001 \times 10^{-3}$
			\\
			\hline
			\hline
			\vspace{0.2in}
		\end{tabular}
	}
	\scalebox{0.6}{
		\begin{tabular}
			%
			%{r@{\hspace*{.30in}}c@{\hspace*{.30in}}c@{\hspace*{.30in}}c@{\hspace*{.30in}}c@{\hspace*{.30in}}c@{\hspace*{.30in}}c} 
			{c c@{\hspace*{.30in}} c@{\hspace*{.30in}} c}
			\multicolumn{4}{c}{Partial fading case II}\\			
			\hline
			\hline
			\# faded history matrices & $\mu =1/10$ & $\mu =1/2$&  $\mu = 9/10$  \\[4pt]
			\hline
			0  & $ 9.26034 \times 10^{-12}$  & $ 2.31391 \times 10^{-11}$  & $ 4.24903 \times 10^{-9}$
			\\
			\hline
			2  & $ 9.27241 \times 10^{-12}$  & $ 2.34361 \times 10^{-11}$  & $ 4.2491 \times 10^{-9}$
			\\
			\hline
			5  & $ 3.37716 \times 10^{-11}$  & $ 6.6476 \times 10^{-10}$  & $ 4.44832 \times 10^{-9}$
			\\
			\hline
			8  & $ 3.8092 \times 10^{-10}$  & $ 1.99961 \times 10^{-8}$  & $ 1.36941 \times 10^{-8}$
			\\
			\hline
			11   & $ 1.47228 \times 10^{-9}$  & $ 2.2715 \times 10^{-7}$  & $ 1.47786 \times 10^{-7}$
			\\
			\hline			
			14   & $ 1.15821 \times 10^{-8}$  & $ 1.64103 \times 10^{-6}$  & $ 1.47098 \times 10^{-6}$
			\\
			\hline
			17   & $ 5.06103 \times 10^{-7}$  & $ 7.87929 \times 10^{-6}$  & $ 1.85274 \times 10^{-5}$
			\\
			\hline
			\hline
			\vspace{0.2in}
		\end{tabular}
	}
	\scalebox{0.6}{
		\begin{tabular}
			%
			%{r@{\hspace*{.30in}}c@{\hspace*{.30in}}c@{\hspace*{.30in}}c@{\hspace*{.30in}}c@{\hspace*{.30in}}c@{\hspace*{.30in}}c} 
			{c c@{\hspace*{.30in}} c@{\hspace*{.30in}} c}
			\multicolumn{4}{c}{Partial fading case III}\\
			\hline
			\hline
			\# faded history matrices & $\mu =1/10$ & $\mu =1/2$&  $\mu = 9/10$  \\[4pt]
			\hline
			0  & $ 9.26034 \times 10^{-12}$  & $ 2.31391 \times 10^{-11}$  & $ 4.24903 \times 10^{-9}$
			\\
			\hline
			2  & $ 9.26023 \times 10^{-12}$  & $ 2.3113 \times 10^{-11}$  & $ 4.24903 \times 10^{-9}$
			\\
			\hline
			5  & $ 1.18462 \times 10^{-11}$  & $ 7.7854 \times 10^{-11}$  & $ 4.23683 \times 10^{-9}$
			\\
			\hline
			8  & $ 1.60656 \times 10^{-10}$  & $ 3.38055 \times 10^{-9}$  & $ 3.65689 \times 10^{-9}$
			\\
			\hline
			11   & $ 1.32413 \times 10^{-9}$  & $ 4.35421 \times 10^{-8}$  & $ 8.84638 \times 10^{-9}$
			\\
			\hline
			14   & $ 7.10271 \times 10^{-9}$  & $ 3.87096 \times 10^{-7}$  & $ 1.7226 \times 10^{-7}$
			\\
			\hline
			17   & $ 2.87023 \times 10^{-8}$  & $ 3.71057 \times 10^{-6}$  & $ 5.12104 \times 10^{-6}$
			\\
			\hline
			\hline
			%			\vspace{0.2in}
			%
		\end{tabular}
	}
\end{table}
\section{\textbf{Summary}}
\label{Sec: Summary}
%%%%%%%%%%%%%%%%%%%%%%%%%%%%%%%%%%%%%%%%%
%

We developed a new $C^{\,0}$-continuous Petrov-Galerkin spectral element method for the problem $\prescript{}{0}{\mathcal{D}}_{x}^{\alpha} u(x) - \lambda u(x) = f(x)$, $\alpha \in (1,2)$, subject to homogeneous boundary conditions. We obtained a weak form, in which the entire fractional derivative load was transferred onto the test functions, allowing us to efficiently employ the standard modal spectral element bases while incorporating Jacobi poly-fractonomials as the test functions. We seamlessly extended the standard procedure of assembling to \textit{non-local assembling} in order to construct the global linear system from local (elemental) mass/stiffness matrices and non-local history matrices. The key to the efficiency of the developed PG method is twofold: i) our formulation allows the construction of elemental mass and stiffness matrices in the standard domain $\left[ -1,1  \right]$ once, and ii) we efficiently obtain the non-local (history) stiffness matrices, in which the non-locality is presented \textit{analytically}. We also investigated local basis/test functions in addition to local basis with global test functions. We demonstrated that the former choice leads to a better-conditioned system and approximability in the spectral element formulation when higher polynomial orders are needed. Moreover, we showed the exponential rate of convergence considering smooth solutions as well as singular solutions with interior singularity; also, the spectral (algebraic) rate of convergence in singular solutions with singularities at boundaries. We also presented the retrieval process of history matrices on uniform grids, which results in faster and more efficient construction and solution of the linear system compared to the on-line computation. In addition, we constructed two non-uniform grids over the computational domain (namely, kernel-driven and geometrically progressive grids), and demonstrated the effectiveness of the non-uniform grids in accurately capturing singular solutions, using fewer number of elements and higher order polynomials. We finally performed a systematic numerical study of non-local effects via both full and partial (history) fading in order to better enhance the computational efficiency of the scheme.

%
%+++++++++++++++++++++++++++++++++++++++
\section*{Acknowledgements}
%+++++++++++++++++++++++++++++++++++++++

The first and the second authors would like to thank Michigan State University for the generous support. The third author was supported by the MURI/ARO on Fractional PDEs for Conservation Laws and Beyond: Theory, Numerics and Applications (W911NF-15-1-0562).

%\newpage
%\bibliographystyle{elsarticle-harv}
%\bibliographystyle{elsarticle-num}
\bibliographystyle{elsarticle-num-names}
\bibliography{RFSLP_Refs2}

\end{document}